\documentclass[hidelinks,onefignum,onetabnum]{siamart220329}
\usepackage[utf8]{inputenc}

\usepackage{amsmath}
\usepackage{amsfonts}
\usepackage{amssymb,latexsym}
\usepackage{mathnotation}
\usepackage{fullpage}
\usepackage[mode=buildnew]{standalone}
\usepackage{tikz}
\usepackage{pgfplots}
\usepackage{standalone}
\usetikzlibrary{math}
\usepackage{xcolor}
\usetikzlibrary{positioning}
\usepackage{catchfile}
\usepackage{tikz-3dplot}
\tdplotsetmaincoords{60}{115}
\pgfplotsset{compat=newest}

\renewtheorem{lemma} {Lemma} [section]
\newtheorem{remark}[lemma]{Remark}
\newtheorem{example}[lemma]{Example}

\def\R{\mathbb{R}}
\def\tg{\cometricg_{\tau}}
\def\bnabla{\dualconnection}
\newcommand{\antidual}{(\solneg*)}

\newcommand{\qdim}{d}

\graphicspath{{TwoLinkPlanarFigure/}}

\def\ctensor{biasing tensor} %

\begin{document}

\title{Optimal control of robotic systems and biased Riemannian splines}
\author{Alejandro Cabrera\thanks{Instituto de Matem\'atica, Universidade Federal do Rio de Janeiro, Rio de Janeiro, Brazil (\email{alejandro@matematica.ufrj.br})} \and Ross L. Hatton\thanks{Robotics and Mechanical Engineering, Oregon State University, Corvallis, OR
(\email{Ross.Hatton@oregonstate.edu})} This work was supported in part by the United States National Science Foundation under awards 1653220 and 1826446 and by the United States Office of Naval Research under award N00014-23-1-2171.}

\maketitle

\begin{abstract}

In this paper, we study mechanical optimal control problems on a given Riemannian manifold $(Q,g)$ in which the cost is defined by a general cometric $\tilde g$. This investigation is motivated by our studies in robotics, in which we observed that the mathematically natural choice of cometric $\tilde g = g^*$---the dual of $g$---does not always capture the true cost of the motion. We then, first, discuss how to encode the system's torque-based actuators configuration into a cometric $\tilde g$. Second, we provide and prove our main theorem, which characterizes the optimal solutions of the problem associated to general triples $(Q,g,\tilde g)$ in terms of a 4th order differential equation.
We also identify a tensor appearing in this equation as the geometric source of "biasing" of the solutions away from ordinary Riemannian splines and geodesics for $(Q,g)$.  Finally, we
provide illustrative examples and practical demonstration of the biased splines as providing the true optimizers in a concrete
robotics system.

\end{abstract}

\setcounter{tocdepth}{2}
\tableofcontents

\section{Introduction}

In mechanical optimal control problems (see~\cite{BulLew}) on a configuration manifold $\configspace$, a natural measure of the cost of following a given trajectory $\gamma(t)\in \configspace$ is the integral of the squared norm of the ``force" $F(t)$ required to follow the trajectory,
\begin{equation} \label{eq:forcenormintro}
C(\gamma) = \int_{t_{0}}^{t_{f}} \| F(t) \|^2_{\cometricg} \ dt.
\end{equation}
In general, calculating this cost involves two metrics on $\configspace$: The first, ``kinematic", metric, $\metricg$, is a Riemannian metric which defines a kinetic-energy norm for configuration velocities $\configdot$ as
\begin{equation} \label{eq:riemannianmetricintro}
\| \configdot \|^2_{\metricg} =  \metricg(\configdot,\configdot) , \ \configdot \in T\configspace .
\end{equation}
and which is additionally used to calculate the both the covariant acceleration $a$ along a trajectory $\gamma$,
\begin{equation} \label{eq:covaccelintro}
a_{\gamma,\metricg} = \nabla^{(g)}_{\configdot} \configdot
\end{equation}
and the associated force covector that must be applied to produce that acceleration,
\begin{equation}  \label{eq:forcedualintro}
F(t) = \metricg(a_{\gamma,\metricg},\bullet).
\end{equation}
(Here, the bullet notation indicates, as usual, an element of the dual vector space obtained by contraction.\footnote{In this case, this operation is often expressed via the ``flat" musical isomporphism, $\metricg^{\flat}(v)\equiv \metricg(v,\bullet) $. Because we later use analogous operations in more general situations, we have chosen to use the bullet notation throughout this text.})
Taking a mechanical linkage as an example (which we expand upon throughout the paper), the Riemannian metric $\metricg$ in~\eqref{eq:riemannianmetricintro} describes the generalized mass for the system, and the covariant acceleration in~\eqref{eq:covaccelintro} describes the extent to which the mass particles in the linkage are accelerating in directions not dictated by the rigid-body and link-hinge constraints. Interesting robotics systems modeled with this setting can be found in~\cite{FraMa,LiFrMa, SmBrFr}.

The second, ``actuation", cometric, which we denote by $\cometricg$, is defined as a quadratic form on covectors which provides the force norm in~\eqref{eq:forcenormintro},
\begin{equation}
\| F \|^2_{\cometricg} = \cometricg(F,F) , \ F \in T^*\configspace .
\end{equation}
Note that $\cometricg$ is independent of the kinematic metric $\metricg$, it encodes information about how and where forces on the system should be measured, and that the selection of $\cometricg$ is an addition of information to the system description.

As discussed below, a classical choice of $\cometricg$ is the  ``dual" $\dualmetricg$ of the system metric $g$, but other situations call for a different choice of $\cometricg$.
Continuing with our mechanical linkage example, the norm of $F$ under the dual of the kinetic-energy Riemannian metric is the root-sum-of-squares of the non-constraint forces felt by the individual particles. If the linkage is actuated from the joints, however, then the cost of generating those forces depends on the configuration-dependent lever-arm distances between the joints and the masses, leading to a different choice of cometric.

\paragraph{Classic dual-metric optimization} Given $\metricg$ on $\configspace$, there is a natural choice for $\cometricg$, namely, to take it as the ``dual" of the first metric, denoted $\dualmetricg$. If we take $\metricg$ and $\cometricg$ as being respectively encoded by metric tensors $\metric$ and $\cometric$, these metric tensors are then related as $\cometric = \dualmetric = \inv{M}$.

With this choice $\cometricg=\dualmetricg$, minimizing the squared norm of the force is equivalent to minimizing the squared norm of its covariant acceleration, %
\begin{equation}
\min_{\gamma} \int_{t_{0}}^{t_{f}} \| F(t) \|^2_{\cometricg} \ dt = \min_{\gamma} \int_{t_{0}}^{t_{f}} \| a_{\gamma,\metricg}(t) \|^2_{\metricg} \ dt
.
\end{equation}
The characteristics of optimal curves in this case are well-understood (see~\cite{CrLe,NoHePa}). When only the initial and end-points are prescribed while the velocities at those points are unconstrained, the least-effort trajectories on the manifold $Q$ are given by $\metricg$-geodesics (shortest, least-curved, and constant-speed). When the initial and final velocities are also prescribed (so that one needs to speed up and slow down, and may need to leave or reach the endpoints with a specified direction), the trajectories are $\metricg$-cubic splines (see~\cite{CrLe,NoHePa,ParkBezier}, and~\cite{BCKS} for a general Hamiltonian approach). In particular, for cases where only the end \emph{speeds} are specified (e.g., the system starts and ends at rest), the $\metricg$-cubic spline follows the same path as the geodesic, only differing from it in the speed profile along the path.

\paragraph{General cometric optimization} In this paper we observe that for many systems, the choice of $\cometricg$ as the dual of $\metricg$ does not accurately model the costs of executing a trajectory, and a different cotangent metric should be used. Under such a cometric $\cometricg$, optimal trajectories may deviate away from the geodesic or cubic spline paths for $\metricg$, either because the $\cometricg$-cost of accelerating through a different region of the configuration manifold is reduced enough to make it worth the extra curvature and pathlength to reach it, or because accelerating obliquely to the geodesic path is cheaper than accelerating directly along it.  It is important to note that these changes in the optimal paths cannot be understood as defining geodesics or splines for some single, possibly different, metric $\metricg'$; rather, they define special curves associated with the pair of metrics $(\metricg,\cometricg)$ through a concrete system of differential equations (see Theorem~\ref{thm:solutions} and Section~\ref{sec:biased} below). We call these solution curves \emph{biased geodesics} and \emph{biased splines} associated with  $(Q,\metricg,\cometricg)$. %

\paragraph{Robot dynamics} The main examples guiding this study come from robot motion planning, in which the model of system dynamics outlined in~\eqref{eq:riemannianmetricintro}--\eqref{eq:forcedualintro} is widely used (albeit typically appearing in the form of the Euler-Lagrange equations of motion, with the connection in the covariant acceleration calculation from~\eqref{eq:covaccelintro} encoded by the Coriolis-centrifugal matrix).  In this context, many planning algorithms assume that least-effort trajectories between given configurations follow the same paths as geodesics for the metric $\metricg$ defined by the system's inertia tensor~\cite{zhang2013kuka,flash1985coordination,biess2007computational,biess2011riemannian}. This assumption is equivalent to taking the effort metric $\cometricg$ as dual to $\metricg$, as above.

In contrast to this common assumption, we present simple examples  of robot motion planning in which properly biased geodesics clearly produce a better optimization than the $\metricg$-geodesic paths. The underlying principle here is that the forces \emph{felt} by the masses in the system (which are measured by $\dualmetricg$) are not necessarily in fixed proportion to the forces \emph{exerted} by the ``actuators" (e.g., motors) which produce the control torques---the effective moment arms from the actuators to the masses depend on both the current configuration of the mechanism and the direction in which force is being applied, and the force felt on each mass particle must be exerted by every actuator between it and the fixed base.

Accordingly, we introduce a cometric $\cometricg$ defined by the way in which the actuators are attached to the system, such that the norm of the force under this cometric captures the actual actuator effort. By means of these choices, the examples then illustrate that planning with biased geodesics and splines for $(Q,\metricg,\cometricg)$ can reduce significantly the cost of controling the system, and that different actuator arrangements can lead to different control strategies for the same mechanism.

\paragraph{Overview of key concepts} \label{sec:overview} The core narrative of the paper is about the contrast between defining the norm for forces on manifolds using the dual of the "mass" metric tensor $g$, as illustrated in Fig.~\ref{fig:spheremetric and dual}, and using a cometric $\tilde g$ that takes into account extra information about the system structure, as illustrated in Fig.~\ref{fig:spheremetricandtorquecometric}.

Our motivating example is mechanical linkages, where, as illustrated in Fig.~\ref{fig:torquesetup}, the ``lever arm'' effects from the links mean that there is not a consistent proportionality between forces felt by the system masses (which is what is measured by the dual metric) and those exerted by the system joints (which we construct our cometric to measure).

When the force norm is defined using the dual metric, splines on the manifold that minimize squared force are equivalent to splines that minimize squared covariant acceleration, and have snap (fourth covariant derivative of position) or tug (second covariant derivative of force) determined by the Riemannian curvature along the trajectory (Equations~\eqref{eq:covariantacceljerksnap},~\eqref{eq:riemannianspline}, and~\eqref{eq:splineforce} below).
The Riemannian curvature term indicates where the divergence/convergence of geodesics on the manifold helps/hurts the spline in reaching its end boundary condition, and therefore whether the introduction of curvature in the spline should be advanced/delayed to (take advantage of)/mitigate the manifold curvature effects.

For non-dual cometrics, the ``force'' covector is replaced by an ``effort" covector which is calculated by dualizing the covariant acceleration using information from the cometric, and the tug equation has a second term that measures the incompatibility between the cometric $\tilde g$ and the manifold’s Riemannian metric $g$ (Equation~\ref{eq:biasedspline}), which guides the spline into regions or along directions where covariant acceleration costs less effort from the system actuation.

Section~\ref{sub:exsbiased} then has several examples of how disconnecting the cometric from the metric in this way affects system trajectories.

\paragraph{Outline of the paper} %
In Section~\ref{sec:kinematics}, we review the geometric setting for the mechanical systems of interest and fix notations. In Section~\ref{sec:torquemetric}, we discuss how to define the metric $\cometricg$ for general actuated mechanical systems and show how to explicitly compute it in concrete examples. In Section~\ref{sec:biased}, we define the optimal control problem associated with the data $(Q,\metricg,\cometricg)$ as above and prove that its solutions come from a Hamiltonian system of equations, following Pontryagin's Maximum Principle (PMP). This system is equivalent to a fourth-order differential equation for the configuration variables and generalizes the equations for ordinary geodesiscs and cubic splines, reducing to them when $\cometricg$ is the dual of $\metricg$. We also work out examples with detail. In Section~\ref{sub:exsbiased}, we describe real-life applications of the general framework in robot motion planning, providing simple examples in which the biased geodesics provide a better optimization than ordinary geodesics.
Some %
useful coordinate expressions are compiled in Appendix~\ref{app:coordinates}.

\paragraph{About the presentation} Although the contents of this paper are mainly mathematical, it is intended to be read by both the geometric-mechanics and the engineering comunities. We thus give detailed procedural accounts of some of the mathematical statements involved, even when some of them are well established. Readers will be, nevertheless, assumed to be familiar with elementary notions of differential geometry. %

\section{Preeliminaries on the geometric setting: kinematics and forces}
\label{sec:kinematics}
In this section, we detail the well-known geometric setting (see e.g. \cite{BulLew}) for the mechanical systems that we will be considering and, at the same time, fix notations. We also include simple explicit examples that will help us illustrate the main constructions of this paper in the later sections. For the interested reader, we also add procedural remarks about how to explicitly compute the ingredients in practical situations.

\subsection{Kinematics: the mass metric}

The general case we consider in this paper is a system whose configuration space is a $\qdim$-dimensional manifold $\configspace$ equipped with a Riemannian metric $\metricg$ that allows for the norm of tangent vectors to be computed as
\begin{subequations}
\begin{align}
\metricg: T_{\config}\configspace \to \euclid, \    \configdot \mapsto \norm{\configdot}_{\metricg}^{2}=\metricg(\configdot,\configdot).
\end{align}
\end{subequations}
Given a choice of local coordinates, this metric can be represented in each tangent space $T_{\config}\configspace$ by a $\qdim\times \qdim$ matrix $\metric(q)$, such that the norm is computed as
\begin{equation}
\norm{\configdot}_{\metricg}^{2}=\metricg(\configdot,\configdot) = \transpose{\configdot} \metric(\config)\, \configdot.
\end{equation}

In many cases of interest,
$\configspace$
is a submanifold of an ambient Riemaniann manifold $\configspace_\amb$ (e.g., $\configspace_\amb=\euclid^N$ with the standard Euclidean metric), $\configspace$ is defined by a set of holonomic constraints, and $\metricg$ is the induced metric from $\configspace_\amb$ onto $\configspace$ (i.e.,
$\norm{\configdot}_g = \norm{J \configdot}_{Q_\amb}$, where $J$ is the Jacobian of the map embedding $\configspace$ into $\configspace_{\amb}$
). %
Moreover, even if $\configspace$ is not inherently a submanifold of a Riemannian manifold, it is always possible in principle to embed it into such an ambient space~\cite{whitney,rodnay2001isometric}.
In this paper, we alternate as is most convenient between the geometries in the ambient space, the coordinate chart, and the manifold embedded into the ambient space,
\begin{equation}
\R^\qdim \supset \varphi(U) \overset{\varphi^{-1}}{\to} \configspace \hookrightarrow \configspace_\amb
\end{equation}
where  $\varphi:U\subset \configspace \to \R^\qdim$ denotes a coordinate chart for $\configspace$.

To fix ideas and notation, let us consider a system composed of $n$ unit point masses in $m$-dimensional space so that $\configspace_\amb = (\R^{m})^n\ni x=(r_1,\dots,r_n), \ r_j \in \R^m$. The total kinetic energy is expressed in terms of the Euclidean metric as
\begin{equation} \label{eq:basicparticleenergy}
\kineticenergy_{\configspace_\amb} = \sum_{i} \frac{1}{2} \transpose{v}_{i} v_{i} = \frac{1}{2} \transpose{v} v,
\end{equation}
where $v_{i} \in \R^m$ is the velocity of the $i$th particle, and $v=(v_1,\dots,v_n) \in T_{x} \configspace_\amb = \R^{m n}$ is the combined velocity vector on the ambient space of all the particles. %
If the point particles are subject to holonomic constraints defining a $\qdim$-dimensional submanifold $\configspace\hookrightarrow \configspace_\amb$ and $v \in T_x \configspace \subset T_x \configspace_\amb$ is a velocity tangent to $\configspace$, the associated kinetic energy is given by
\begin{equation}
\kineticenergy_{\configspace} = \frac{1}{2} \metricg(v,v) = \frac{1}{2}v^T v,
\end{equation}
where the second equality follows from the fact that the metric $\metricg$ on $\configspace$ is induced from that of the ambient space.
We can also choose a $\qdim$-dimensional chart $\varphi:U\subset \configspace\to  \R^\qdim$ so that the configuration is parameterized by $q=(q^1,\dots,q^\qdim)\in \R^\qdim$ via
\begin{equation}
(\R^m)^n \ni x = \varphi^{-1}(q)=(\varphi^{-1}_1(q),\dots,\varphi^{-1}_n(q)); \ \ r_i=\varphi^{-1}_i(q)\in \R^m, \ i=1,\dots, n.
\end{equation}
The corresponding coordinate expression for the induced kinetic energy function on $\configspace$ can be obtained as follows: The velocity vectors in $\configspace_\amb$ are given by the associated Jacobian matrices via
\begin{equation} \label{eq:defJac}
v_i = \overset{J_i\equiv J_i(q)}{\overbrace{(\partial \varphi^{-1}_i/ \partial q^j)(q)}} \ \configdot^j, \ \configdot\in \R^\qdim, \ i=1,\dots, n. \text{ (Using Einstein's summation convention.)}
\end{equation}
With these notations, the induced kinetic energy function in the chart space $\varphi(U)$ is given by
\begin{align}
\kineticenergy_{\varphi(U)}    &= \frac{1}{2} \transpose{\configdot} \underbrace{\left(\sum_{i} \transpose{\jac_{i}(q)} \jac_{i}(q)\right)}_{M(q)} \configdot \label{eq:kinenQsetting},
\end{align}
where the matrix $J(q):\R^\qdim \to (\R^m)^n$ is seen as a linear transformation from vectors expressed in the tangent space of $Q$ to their corresponding vectors in the tangent space of the ambient space $Q_\amb$.

One can interpret $M\equiv M(q)$ physically as the effective mass matrix of the particles with respect to the configuration variables $q$. We also recall that this pullback procedure can be generalized to non-unit point masses by introducing a weighting term into the summation in~\eqref{eq:basicparticleenergy}, and to systems with continuous mass distributions by replacing the summations in the above equations with integrals.

\begin{remark}\textsc{(Other ambient spaces $\configspace_\amb$)}
There is an important class of applications (e.g., when modeling a robot constructed of discrete links) in which is better to consider the above setting with $\configspace_\amb$ given by $SE(2)^{k}$ or $SE(3)^{k}$ and representing configurations of $k$ rigid bodies. In these cases, $\configspace\hookrightarrow \configspace_\amb$ is defined by the ``interlink" holonomic constraints  (e.g., pivoting or sliding joints) that restrict the relative positions and orientations of these rigid bodies to a $\qdim$-dimensional manifold. The metric on $\configspace$ is then defined by the pullback of a natural metric on $\configspace_\amb$, which encodes the  classical mass and moment of inertia tensors of these bodies, and is itself defined by the pullback of the Euclidean metric on the particle motions through the rigid body constraints.
\end{remark}

We now begin with a series of examples involving the sphere $S^2\hookrightarrow \R^3$, which will help us illustrate each step towards our biased geodesics and splines.

\begin{example}\textsc{(Mass matrix on a sphere)}\label{ex:massmatrixsphere}
Let us consider the sphere of radius 1 in 3-dimensional space,
\begin{equation}
\configspace=\sphere^2:=\{ (x,y,z): x^2+y^2+z^2=1\}
\hookrightarrow \configspace_\amb = \R^3.
\end{equation}
Physically, one can think of a unit mass point particle is subject to holonomic constraints restricting its motion to the surface $\configspace$.
We now detail how to obtain the pullback metric $\metricg$ from the Euclidean metric on the ambient space.  Let us consider a coordinate chart $\varphi$ as above given by the longitude-latitude coordinates $q=(\lon,\lat)$,
\begin{equation}
(x,y,z) = \varphi^{-1}(\lon,\phi) := (\cos{\lon}\cos{\lat},\ \sin{\lon}\cos{\lat},\ \sin{\lat}).
\end{equation}
The associated Jacobian matrix $J(q):\R^2\to \R^3$ is given by
\begin{equation}
J(q)=
\begin{bmatrix}
-\cos{\lat}\sin{\lon} &  -\cos\lon\sin\lat \\
\phantom{-}\cos{\lat}\cos{\lon} & -\sin\lat\sin\lon \\
\phantom{-}0 & \phantom{-}\cos{\lat}
\end{bmatrix}
\end{equation}
and thus the associated mass matrix $\metric_{\sphere}(q):\R^2 \to \R^2$ yields
\begin{figure}
\centering
\includegraphics[width=\textwidth]{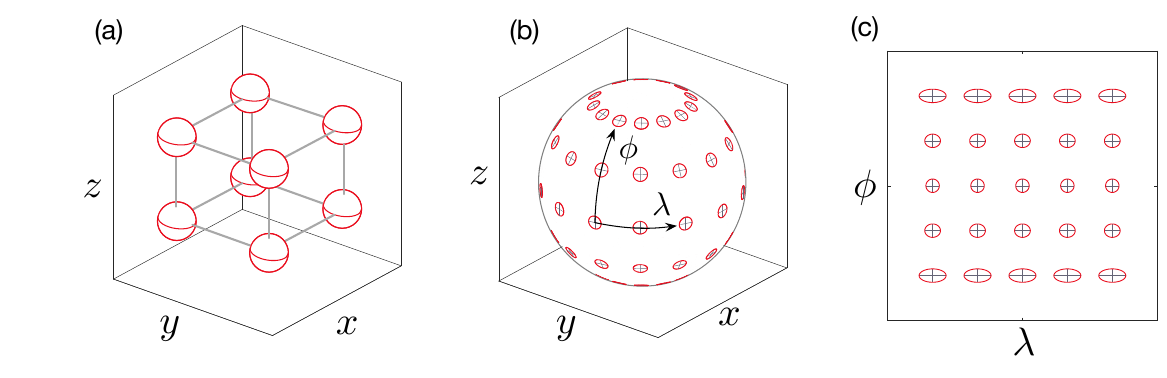}
\caption{%
(a) The Tissot indicatrices for the Euclidean metric in $\R^3$ are spheres, indicating equal weighting for all coordinate components of a vector when calculating its length. Connecting lines between spheres are a coordinate grid to help with 3D viewing perspective. (b) Restricting the Euclidean metric to the tangent spaces on a sphere produces a set of circles. (c) The indicatrices for the induced metric on the longitude-latitude coordinate chart is a set of longitudinally stretched ellipses, as if the circles on the sphere had been stretch-flattened onto the plane.}
\label{fig:metricholonomicrestriction}
\end{figure}
\begin{equation}
\metric_{\sphere}(q)=J(q)^T J(q) = \begin{bmatrix} \cos^{2}\lat & 0 \\ 0 & 1 \end{bmatrix} \label{eq:spheremetric},
\end{equation}
which determines the kinetic energy on $\configspace$ as seen in the chart $\varphi$ via eq.~\eqref{eq:kinenQsetting}. For later use, in Fig.~\ref{fig:metricholonomicrestriction} we include a graphical representation of the metric defined by $M(q)$ on this chart using generalized Tissot indicatrices.

A given change in the longitude $\lon$ results in a relatively smaller motion in the ambient space when the latitude $\lat$ is large (near the poles) than when it is small (near the equator). Conversely, a given motion on the surface of the sphere requires more change in $\lon$ when $\lat$ is large than when it is small.
This property can be visualized in terms of the Tissot indicatrices (sets of metric-unit-length velocities), which at high latitude are ``stretched" in the longitudinal direction, as illustrated in Fig.~\ref{fig:metricholonomicrestriction}(c).%
\end{example}

\subsection{Acceleration}

Given a trajectory $\gamma(t)\in Q$, the \emph{covariant acceleration} $a_{\gamma}$ describes the ``true" rate at which the velocity $\dot{\gamma}$ changes along the curve. In the case $Q\hookrightarrow Q_\amb$ above, calculating the covariant acceleration corresponds to taking the second derivative in $Q_\amb$ and then projecting to the tangent space to $Q$. The projection operation annihilates the acceleration normal to the embedded manifold so that $a_{\gamma}$ can be considered an intrinsic representation of the extent to which $\gamma$ is changing speed or turning within $Q$.

\smallskip

\noindent {\bf Calculating covariant acceleration intrinsically.}
In general, $a_{\gamma}$ can be calculated directly in $Q$-coordinates (without actually constructing an isometric embedding) as the covariant derivative of the trajectory velocity, i.e. the derivative with respect to the Levi-Civita connection $\nabla\equiv \nabla^g$ on $(Q,g)$,
\begin{equation}\label{eq:defcova}
a_{\gamma,\metricg}(t) = \nabla_{\dot{\gamma}(t)}\dot{\gamma}(t).
\end{equation}
As reviewed in Appendix~\ref{app:coordinates}, this operation combines the standard coordinate acceleration $\ddot{\config} = \frac{d^2}{dt^2}\gamma(t)$ with a correcting term computed from the metric tensor and its derivatives. This correcting term can be thought of as compensating for any ``false" acceleration contributed to $\ddot{\config}$ by the ways in which the coordinate grid is stretched or curved as it is overlaid on the manifold $Q$.

\begin{example}\textsc{(Acceleration on the sphere)}
In the setting of Example~\ref{ex:massmatrixsphere}, let $P_x: \R^3 \to T_x \sphere^2$ be the orthogonal projection onto the tangent space at $x\in \sphere^2$ (seen naturally as a vector subspace in $\R^3$). A point moving over the surface of the sphere along a curve $\gamma(t)\in \sphere^2\hookrightarrow \R^3$ has covariant acceleration given by
\begin{equation}
a_{\gamma,\metricg} (t) = P_{\gamma(t)}(\ddot \gamma (t)) \in T_{\gamma(t)} \sphere^2.
\end{equation}
We now consider some special cases, for later use.
If the trajectory $\gamma(t)$ is a geodesic in $(\sphere^2,\metricg)$, namely it follows great circle at constant speed (as in Fig.~\ref{fig:sphereacceleration}(a)), then the ambient acceleration is entirely normal to the surface of the sphere, and hence $a_\gamma(t) = 0$.
\begin{figure}
\centering
    \begin{tikzpicture}[tdplot_main_coords, scale = 2.5]
\coordinate (P) at ({1/sqrt(3)},{1/sqrt(3)},{1/sqrt(3)});
 
 \draw[dashed, gray] (0,0,0) -- (-1,0,0);
\draw[dashed, gray] (0,0,0) -- (0,-1,0);
 
\draw[-stealth] (0,0,0) -- (1.80,0,0) node[below left] {$x$};
\draw[-stealth] (0,0,0) -- (0,1.30,0) node[below right] {$y$};
\draw[-stealth] (0,0,0) -- (0,0,1.30) node[above] {$z$};
 
\shade[ball color = lightgray, opacity = 0.5] (0,0,0) circle (1cm);
 
%\tdplotsetrotatedcoords{0}{0}{0};
%\draw[dashed, tdplot_rotated_coords, gray] (0,0,0) circle (1);
 
%\tdplotsetrotatedcoords{0}{0}{0};
%\draw[thick,tdplot_rotated_coords, red] (0,0,{1/sqrt(2)}) circle ({1/sqrt(2)});

\tdplotsetrotatedcoords{-230}{45}{0};
\draw[thick,tdplot_rotated_coords, red] (0,0,0) circle (1);

\draw[-stealth,thick,tdplot_rotated_coords] (-1,0,0) -- (-1,-.75,0) node[left] {$v(t)$};
%\draw[-stealth,thick,tdplot_rotated_coords] (-1,0,0) -- (-1,0,.25) node[right] {$a(t)$};
\draw[fill = red,tdplot_rotated_coords] (-1,0,0) circle (0.5pt);

\tdplotsetrotatedcoords{-230}{45}{30};
\draw[-stealth,thick,tdplot_rotated_coords] (-1,0,0) -- (-1,-.65,0)  node[ left ] {$v(t+\Delta t)$ };
\draw[fill = lightgray!50,tdplot_rotated_coords] (-1,0,0) circle (0.5pt);

%\tdplotsetrotatedcoords{-60}{45}{0};
%\draw[dashed,tdplot_rotated_coords, black] (0,0,0) circle (1);

%\tdplotsetrotatedcoords{90}{90}{90};
%\draw[dashed, tdplot_rotated_coords, gray] (1,0,0) arc (0:180:1);
% 
%\tdplotsetrotatedcoords{0}{90}{90};
%\draw[dashed, tdplot_rotated_coords, gray] (1,0,0) arc (0:180:1);
% 

%\draw[thick, -stealth] (0,0,0) -- (P) node[right] {$P$};
% 
%\draw[thin, dashed] (P) --++ (0,0,{-1/sqrt(3)});
%\draw[thin, dashed] ({1/sqrt(3)},{1/sqrt(3)},0) --++
%(0,{-1/sqrt(3)},0);
%\draw[thin, dashed] ({1/sqrt(3)},{1/sqrt(3)},0) --++
%({-1/sqrt(3)},0,0);
% 
%\draw[fill = lightgray!50] (P) circle (0.5pt);
\end{tikzpicture}
    \begin{tikzpicture}[tdplot_main_coords, scale = 2.5]
\coordinate (P) at ({1/sqrt(3)},{1/sqrt(3)},{1/sqrt(3)});
 
 \draw[dashed, gray] (0,0,0) -- (-1,0,0);
\draw[dashed, gray] (0,0,0) -- (0,-1,0);
 
\draw[-stealth] (0,0,0) -- (1.80,0,0) node[below left] {$x$};
\draw[-stealth] (0,0,0) -- (0,1.30,0) node[below right] {$y$};
\draw[-stealth] (0,0,0) -- (0,0,1.30) node[above] {$z$};
 
\shade[ball color = lightgray, opacity = 0.5] (0,0,0) circle (1cm);
 
%\tdplotsetrotatedcoords{0}{0}{0};
%\draw[dashed, tdplot_rotated_coords, gray] (0,0,0) circle (1);
 
\tdplotsetrotatedcoords{0}{0}{0};
\draw[thick,tdplot_rotated_coords, red] (0,0,{1/sqrt(2)}) circle ({1/sqrt(2)});

\tdplotsetrotatedcoords{-230}{45}{0};
\draw[dashed,tdplot_rotated_coords, red] (0,0,0) circle (1);

\draw[-stealth,thick,tdplot_rotated_coords] (-1,0,0) -- (-1,-.75,0) node[below] {$v(t)$};
\draw[-stealth,thick,tdplot_rotated_coords] (-1,0,0) -- (-1,0,.25) node[right] {$a(t)$};
\draw[fill = red,tdplot_rotated_coords] (-1,0,0) circle (0.5pt);

\tdplotsetrotatedcoords{-205}{0}{0};
\draw[-stealth,thick,tdplot_rotated_coords] (-{1/sqrt(2)},0,{1/sqrt(2)}) -- (-{1/sqrt(2)},-.65,{1/sqrt(2)}) node[ right, ] {$v(t+\Delta t)$ };
\draw[fill = lightgray!50,tdplot_rotated_coords] (-{1/sqrt(2)},0,{1/sqrt(2)}) circle (0.5pt);

%\tdplotsetrotatedcoords{-60}{45}{0};
%\draw[dashed,tdplot_rotated_coords, black] (0,0,0) circle (1);

%\tdplotsetrotatedcoords{90}{90}{90};
%\draw[dashed, tdplot_rotated_coords, gray] (1,0,0) arc (0:180:1);
% 
%\tdplotsetrotatedcoords{0}{90}{90};
%\draw[dashed, tdplot_rotated_coords, gray] (1,0,0) arc (0:180:1);
% 

%\draw[thick, -stealth] (0,0,0) -- (P) node[right] {$P$};
% 
%\draw[thin, dashed] (P) --++ (0,0,{-1/sqrt(3)});
%\draw[thin, dashed] ({1/sqrt(3)},{1/sqrt(3)},0) --++
%(0,{-1/sqrt(3)},0);
%\draw[thin, dashed] ({1/sqrt(3)},{1/sqrt(3)},0) --++
%({-1/sqrt(3)},0,0);
% 
%\draw[fill = lightgray!50] (P) circle (0.5pt);
\end{tikzpicture}
 

\caption{Draft of figure illustrating covariant acceleration on a sphere. Will also include what these look like on the coordinate chart. If you've got experience with drawing nicer spheres in tikz, I'd like to learn.}
\label{fig:sphereacceleration}
\end{figure}
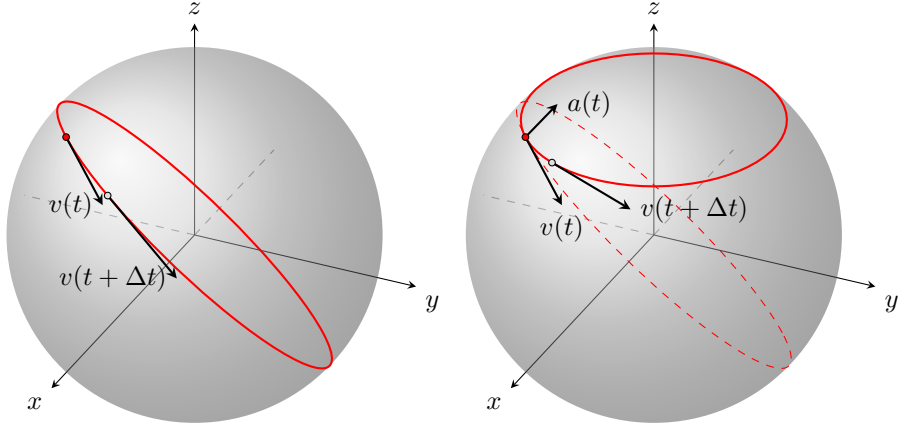
If it was changing speed, then $a_\gamma(t)$ would be a non-vanishing vector parallel to the velocity $\dot \gamma (t)$. Along a non-great-circle path, e.g., following the $\lat = \pi/4$ latitude line at constant speed as in Fig.~\ref{fig:sphereacceleration}(b), the covariant acceleration $a_\gamma$ captures the acceleration required to ``steer" the system off of the great circle to which it is currently tangent.
This acceleration can be found in longitude-latitude coordinates by %
considering $\config(t)=(\lon(t),\pi/4)$ with the metric from~\eqref{eq:spheremetric}, yielding
\begin{equation}
a_{\gamma,\metricg} = \begin{bmatrix}
0\\
\dot{\lon}^{2}/2
\end{bmatrix}.
\end{equation}
\end{example}

\smallskip

\noindent {\bf Magnitude of covariant acceleration.}
Calculating the magnitude of the covariang acceleration as its $\metricg$-norm,
\begin{equation}
\norm{a_{\gamma,\metricg}}^{2} = \metricg(a_{\gamma,\metricg},a_{\gamma,\metricg}),
\end{equation}
coincides with the ambient magnitude of the tangential component of the acceleration of the isometrically embedded trajectory with respect to $Q\hookrightarrow Q_\amb$. For systems composed of point masses or rigid bodies, whose metrics are constructed as in~\eqref{eq:kinenQsetting}, $\norm{a_{\gamma,\metricg}}^{2}$ is additionally equal to the mass-weighted sum of the squared accelerations of the individual particles.

\smallskip

\noindent {\bf Higher-order derivatives.}
Geometrically meaningful forms of the higher-order derivatives of position (e.g., jerk and snap) can be calculated by repeatedly applying the covariant derivative operator. For example, the covariant jerk is defined by,
\begin{equation}
j_{\gamma,\metricg}(t) = \nabla_{\dot{\gamma}(t)}a_{\gamma,\metricg}(t).
\end{equation}
Intuitively, the covariant jerk is zero when the magnitude of the covariant acceleration and its angle with respect to the trajectory (both measured either in isometric ambient coordinates or using $\metricg$ in working coordinates) remains constant (i.e., $a$ is parallel-transported along $\gamma$).

\subsection{Forces}
The customary generalization of Newton's second law, $F=ma$, consists in identifying the force associated with a trajectory $\gamma(t)\in \configspace$ as being dual with respect to $\metricg$ of the covariant acceleration,
\begin{equation}\label{eq:forcedualaccel}
F_{\gamma,\metricg} = \metricg(a_{\gamma,\metricg},\bullet) \in T^*_{\gamma(t)}\configspace.
\end{equation}
(Note that the force associated with a trajectory $\gamma$ is thus a field of covectors along $\gamma$.)
Given a coordinate chart $\varphi$ and the associated coordinate trajectory $q(t)=\varphi(\gamma(t))\in \varphi(U)$, the associated force on $\varphi(U)$ is then given by
\begin{equation} \label{eq:intrinsicforceaccel}
F_{\gamma,\metricg}(t) = M(q(t))\, a_{\gamma,\metricg}(t)
\end{equation}
where $M(q)$ is the metric tensor at the current configuration, and $a_{\gamma,\metricg}$ is the covariant acceleration computed via eq.~\eqref{eq:defcova}.

\begin{example}\textsc{(Forces on the sphere)}\label{ex:forcessphere}
Here we continue with the setting of Example~\ref{ex:massmatrixsphere}. The force required to produce a covariant acceleration $a_{\gamma,\metricg}=( a^{\lon},  a^{\lat})$ on a sphere is
\begin{equation}
F_{\gamma,\metricg} =
\begin{bmatrix}
F_{\lon} \\ F_{\lat}
\end{bmatrix}
=
\overset{\metric_{\sphere}(q)}{\overbrace{\begin{bmatrix} \cos^{2}\lat & 0 \\ 0 & 1 \end{bmatrix}}}
\begin{bmatrix}
a_{\lon} \\ a^{\lat}
\end{bmatrix}
=
\begin{bmatrix}
(\cos^{2}{\lat}) a^{\lon} \\ a^{\lat}
\end{bmatrix}.
\end{equation}
Observe that the $cos^{2}{\lat}$ term in $\metric$ means that for a fixed value of $a^{\lon}$, the corresponding force component $F_{\lon}$ is small when $q=(\lon,\lat)$ represents a point near the poles. This scaling corresponds to the length scale associated with $\lon$ being proportional to $1/\cos{\lat}$, which means that as $\lat$ increases, a given $a^{\lon}$ represents a smaller ``real" acceleration in the ambient space,  $F_{\lon}$ represents a \emph{larger} ``real" force tangent to the sphere, and that these scaling factors compound with each other rather than canceling out. %

\end{example}

\smallskip

\noindent {\bf Derivatives of forces.}\label{sec:forcederivative}
The true rate at which the force changes across a trajectory---again, correcting for any ``false" changes introduced by distortions of the coordinates as they are overlaid on the manifold---can be found using the dual covariant derivative on the manifold. For example, the ``yank" (first derivative of the force with respect to time~\cite{LinYank}) is calculated as
\begin{equation}
y_{\gamma,\metricg}(t) = \dualconnection_{\dot{\gamma}(t)}F_{\gamma,\metricg}(t),
\end{equation}
where $\nabla^*$ denotes the linear connection on $T^*Q$ which is dual to the given Levi-Civita $\nabla$ on $TQ$. %
One way of defining $\nabla^*$ is via the equation
\begin{equation}
L_{X}(\covprod{\omega}{v}) = \covprod{\dualconnection_{X}\omega}{v} +  \covprod{\omega}{\nabla_{X}v}, \qquad X,v\in \mathfrak{X}(Q), \omega \in \Omega^1(Q),
\end{equation}
so that the directional derivative of the product of the two fields along any vector field $X$ reflects only the true changes in the fields.
As reviewed in Appendix~\ref{app:coordinates}, the correcting terms in the dual covariant derivative take a slightly different form than they do in the standard covariant derivative. These differences correspond to the inverse ways in which coordinate grid distortion affects the coordinate representations of vectors and covectors.

\smallskip

\noindent {\bf Cometrics and dual metrics.}
To calculate the magnitude of the force covectors, we must introduce a second metric structure, a \emph{cometric} for the cotangent bundle $T^*\configspace$, denoted $\cometricg$. This cometric provides a norm for cotangent vectors such as forces,
\begin{align}
\cometricg: T_{x}^{*}\configspace \to \euclid, \
F \mapsto \norm{F}_{\cometricg}^{2}=\cometricg(F,F),
\end{align}
and is encoded in coordinates by a matrix $\cometric(\config)$ such that
\begin{equation} \label{eq:dualnormconsistency}
\norm{F}_{\cometricg}^{2}=\cometricg(F,F) = F \cometric(\config)\, \transpose{F},
\end{equation}
using the convention that covectors are written as rows.

The canonical choice for $\cometricg$ is the \emph{dual metric} $\dualmetricg$ to $\metricg$, defined such that applying the dual metric to the $\metricg$-dual of a vector is equivalent to applying $\metricg$ directly to the vector,
\begin{equation} \label{eq:defgdual}
\dualmetricg\bigl(\metricg(v,\bullet),\metricg(v,\bullet)\bigr) = \metricg(v,v), \qquad \ v \in T_x \configspace.
\end{equation}
The matrix encoding $\dualmetricg$ is thus the inverse $M^{-1}$ of the matrix $M$ encoding $\metricg$: Taking $F = \metricg(a,\bullet)$ and $\metric$ known to be symmetric, we have
\begin{equation}
\norm{F}_{\dualmetricg}^{2} = F\inv{\metric}\, \transpose{F}
= \transpose{(\metric a)} \inv{\metric} (\metric a)
= \transpose{a} (\transpose{\metric} \inv{\metric} \metric) a
= \transpose{a} \metric a
= \norm{a}_{\metricg}^2.
\end{equation}
As with the metric $\metricg$, the dual metric $\dualmetricg$ can be obtained by pullback from the ambient space along $\configspace\hookrightarrow \configspace_\amb$, and interpreted as being the restriction of the (Euclidean) cometric on the ambient space onto the manifold. %
In coordinates, we can construct this pullback by noting that the map from a covector expressed in coordinate bases to its corresponding minimal (i.e., tangent-projected) ambient covector is given by the pseudoinverse $J^+$ of the Jacobian $J$ of the embedding function (see eq.~\eqref{eq:defJac}).
\begin{equation}
\norm{F}_{\dualmetricg}^{2} = F_{\text{amb}}\ \transpose{F_{\text{amb}}}
= (F \pinv{\jac}) \transpose{(F\pinv{\jac})}
= F \underbrace{(\pinv{\jac} \pinvtranspose{\jac})}_{\displaystyle{\inv{\transpose{\jac}\jac}}} \transpose{F}
=F\inv{\metric} \, \transpose{F}.
\end{equation}
Finally, the dual connection $\nabla^*$ and the dual metric $\dualmetricg$ are related by the compatibility relation
\begin{equation} \label{eq:compatibility}
L_{X} \dualmetricg(\onef_{1},\onef_{2}) = \dualmetricg(\dualconnection_{X}\onef_{1},\onef_{2}) + \dualmetricg(\onef_{1},\dualconnection_{X}\onef_{2}), \qquad X\in\mathfrak{X}(Q), \omega_1,\omega_2 \in \Omega^1(Q).
\end{equation}
This relationship follows directly from the compatibility between $\nabla$ and $g$ by duality: The dual metric is ``unchanging" with respect to the metric, so any change in the dual-metric product across the space depends only on ``real" changes in the two input covector fields.

\begin{figure}
\centering
\includegraphics[width=\textwidth]{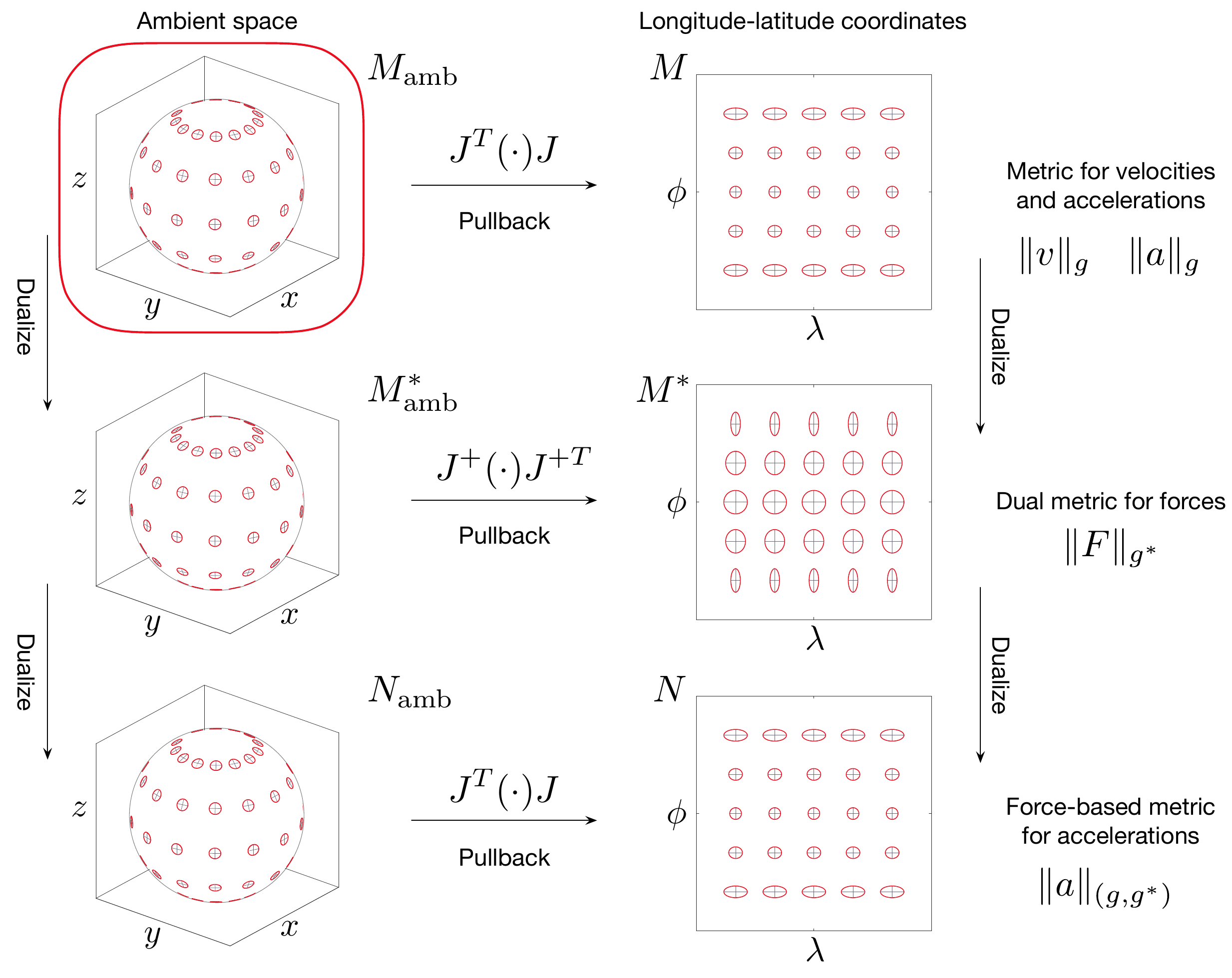}
\caption{The metric, dual metric, and force-based acceleration metric for the sphere are all generated via pullback and dualization operations stemming from the restriction of the ambient-space Euclidean metric onto the sphere. The norm $\norm{a}_{(\metricg,\dualmetricg)}$ and its metric tensor $\inmetric$ are defined in eq.~\eqref{eq:defhfromtg} and~\eqref{eq:defhfromtgN}. The dualization operations are defined by the standard Euclidean ``identity" dualizations on the ambient space, and can be constructed on the coordinate chart by composition of the ambient-space dualizations with the relevant pullback operations.}
\label{fig:spheremetric and dual}
\end{figure}
\begin{example}\textsc{(Force metrics on the sphere)}
Continuing with Example~\ref{ex:forcessphere}, the dual metric $\dualmetricg$ on the sphere takes the ambient basis as orthonormal (i.e., having an identity cometric tensor), and, as illustrated in Fig.~\ref{fig:spheremetric and dual}, pulls back into longitude-latitude coordinate space as
\begin{align}
\dualmetric_{\sphere}(q) = \inv{\metric}_{\sphere}(q) =
\begin{bmatrix}
\frac{1}{\cos^{2}\lat} & 0 \\ 0 & 1
\end{bmatrix}
.\nonumber
\end{align}
The $1/\cos^{2}{\lat}$ term corresponds to the property that (as discussed in Example~\ref{ex:forcessphere}) the ambient length scale associated with $\lon$ becomes smaller as $\lat$ increases, and so the ``real" force associated with a given value of $F_{\lon}$ increases proportionally (and the amount of $F_{\lon}$ for a given ``real" force decreases proportionally).  In the illustration of $\dualmetric_{\sphere}$ at the middle-right of Fig.~\ref{fig:spheremetric and dual}, the ellipses of unit force are consequently ``narrow" at large $\lat$ (because in this region it takes little $F_{\lon}$ [torque on the polar axis] to produce a unit force [as measured in the ambient space at the surface of the sphere]).
The on-manifold Tissot indicatrices in the left column of Fig.~\ref{fig:spheremetric and dual} illustrate the property that the metric, co-metric, and inverse-dualized cometric are all orthonormal on the manifold, such that the norms of velocity, acceleration, force, and acceleration-in-terms-of-force encoded by these metric structures are all equal to their values as would be calculated in the ambient space.

\end{example}

\section{The torque-cost cometric of an actuator-controlled system}\label{sec:torquemetric}

In this section, we explore the idea that it can be crucial to use a cometric $\cometricg$ that is \emph{not} the dual of the manifold metric $\metricg$. We first introduce a simple geometric description of what it means to choose a different cometric, along with notation for working with metrics and cometrics that are not dual to each other. We then discuss specific applications in which a non-dual cometric naturally appears: modeling mechanical systems in which forces are supplied by \emph{actuators} attached to the physical sructure at specific points. This leads us to a torque-cost cometric $\tg$. We discuss how to characterize this cometric in practical situations. Finally, we present a set of examples illustrating how this torque-based cometric  biases the cost of generating system accelerations and discuss practical methods involving the computation of $\tg$.

\subsection{General cometrics and their induced metrics}

The choice of the dual metric $\dualmetricg$ for the cometric $\cometricg$ is conventional (to the point where many sources take ``dual metric" and ``cometric" as synonyms), but it is not mandatory. Using $\dualmetricg$ as the relevant cometric is actually an active declaration that the relevant magnitude of a force is its (Euclidean) length in a (dual) basis orthonormal with respect to the manifold geometry encoded in $\metricg$. If the manifold is defined within an ambient space, an equivalent description of $\dualmetricg$ is that it takes the magnitude of the force as the length of the ambient-space force component which is cotangent to the manifold.

The key point is that if we have access to more informaton about the system structure (e.g., the ``leverage" that the actuators have on a mechanical system in different configurations, as discussed in \S\ref{sub:tg} below), we can use this information to construct a cometric $\cometricg$ (with metric tensor $\cometric$) that better reflects the costs of system motions.

We also describe in this subsection two associated quantities. The first quantity is an associated ``induced metric" $\inmetricg$ (with metric tensor $\inmetric$) which provides a norm of accelerations corresponding to the cost of producing them, i.e., such that $\inmetricg(a,a) = \cometricg(F,F)$. The second quantity is an ``effort" covector $E$ produced by dualizing the covariant acceleration via $\inmetricg$; this covector has the property that its product with the covariant acceleration is equal to the cometric norm of the force, $\covprod{E}{a} = \cometricg(F,F)$. In \S\ref{subsec:biasedsplines} below, we use this notion of effort to understand the form of optimal trajectories through the configuration space.

\subsubsection{Constructing the induced metric $h$}
Once a cometric $\cometricg$ has been chosen, the norms of forces can be taken with respect to $\cometricg$ as in eq.~\eqref{eq:dualnormconsistency}. Combining this force-norm with the relation $F=g(a,\bullet)$ between force and acceleration produces a $(\metricg,\cometricg)$-induced  acceleration-norm $\norm{\bullet }_{(\metricg,\cometricg)}$ that describes the force-cost of producing a given acceleration,
\begin{equation}\label{eq:defhfromtg}
\norm{a}_{(\metricg,\cometricg)} := \cometricg(F,F)
=\cometricg\bigl(\metricg(a,\bullet),\metricg(a,\bullet)\bigr)
=: \inmetricg(a,a),
\end{equation}
where the last equality is the definition of the underlying induced metric $\inmetricg$. In coordinates, the induced metric evaluates as
\begin{equation} \label{eq:defhfromtgN}
\norm{a}_{(\metricg,\cometricg)} = \transpose{(\metric a)} \cometric (\metric a)
=\transpose{a} (\metric \cometric \metric) a
=\transpose{a} \inmetric a,
\end{equation}
and we call the matrix $\inmetric=M\cometric M$ that encodes $\inmetricg$ the \emph{anti-dual of $\cometric$ with respect to $\metric$}.

Note that in the case that the cometric is the dual metric, ~\eqref{eq:defgdual} gives the induced metric as being equal to the original metric,
\begin{equation}\label{eq:hequalsg}
\cometricg=\dualmetricg \Rightarrow h=g,
\end{equation}
and in coordinates,
\begin{equation}
\cometric=\inv{\metric} \Rightarrow N=M.
\end{equation}

\subsubsection{Effort covector} \label{sec:effortcovectordescrip}

We define the effort required to produce a covariant acceleration $a$ as the covector $E$ whose product with the acceleration is equal to the cost of producing that acceleration,
\begin{equation}
\covprod{E}{a} = \cometricg(F,F).
\end{equation}
It follows that this effort $E$ is the dual of the acceleration with respect to the induced metric $h$,
\begin{equation}
E = \inmetricg(a,\bullet).
\end{equation}
The magnitude of the effort can be calculated with respect to the cometric $h^*$ which is dual to $h$ so that
\begin{equation}
h^*(E,E) = \cometricg(F,F).
\end{equation}
This equality can be checked using the coordinate relations $E=Na$, $N=M\cometric M$, $F=Ma$ and taking the matrix $N^*$ corresponding to $h^*$ as $N^{*}=N^{-1}$.

Note that when $\cometricg=\dualmetricg$ is chosen as the dual metric, relation~\eqref{eq:hequalsg} implies that the effort is equal to the standard force, $E=F$. In this case, $h^*=g^*$ also reduces to the standard cometric.

\medskip

\noindent {\bf Immediate implications of general cometrics on optimal trajectories.} In Section~\ref{sec:biased}, we will introduced an optimal control problem associated with the data $(Q,g,\tilde g)$.
The property that $\inmetric\neq M$ and $E \neq F$ for cometrics $\cometricg$ that are not dual to $\metricg$ forms the crux of our analysis in this paper: it means that accelerations which are of equal length under the standard norm $\norm{a}_{\metricg}$ can have different force costs $\norm{a}_{(\metricg,\cometricg)}=h(a,a)$ depending on where in the configuration space they are located, or on the direction in which they are pointing. This will be the source of the ``bias" that the solutions to the control problem present with respect to the standard control problems defined by taking $\tilde g = g^*$.

Additionally, as in several places in this manuscript, we emphasize that the optimal control problem resulting from the introduction of $\cometricg$ to measure costs is \emph{not} equivalent to a standard optimal control problem in which the underlying metric is simply replaced with the induced metric $\inmetricg$. The true geometric description is given by the triple $(\configspace,\metricg,\cometricg)$ and the information of each of the two metrics is needed separately as they provide independent pieces of information about the system: $\metricg$ providing the kinematics which allow for the calculation of accelerations and relate accelerations to forces, and $\cometricg$ providing an accurate cost of producing each force due to the system's construction. The induced metric structure $\inmetricg$ is, however, useful for some operations in defining the optimal curves, as we describe in \S\ref{sec:biased}.

\subsubsection{Metric compatibility of general cometrics}

In eq.~\eqref{eq:compatibility}, we noted that the dual metric $\dualmetricg$ is compatible with the geometry defined by $(Q,g)$ through the dual connection $\nabla^*$. General cometrics $\cometricg$ (including induced cometrics $h^*$) are not compatible with the geometry defined by $(Q,\metricg)$, and the Lie derivative of $\cometricg(\onef_{1},\onef_{2})$  additionally includes a term which captures how the cometric changes in relation to $(Q,g)$,
\begin{equation}\label{eq:incompatibility} L_{X}\cometricg(\onef_1,\onef_2) = \overbrace{(\dualconnection \cometricg)_X(\onef_1,\onef_2)}^{\mathclap{\text{change in cometric relative to metric}}} + \cometricg(\dualconnection_X \onef_1, \onef_2) + \cometricg( \onef_1, \dualconnection_X \onef_2). \end{equation}
This change in the cometric across the manifold plays a key role in our analysis later in this paper.

Equivalently, eq. ~\eqref{eq:incompatibility} can be seen as the definition of a ``\ctensor" between the cometric $\cometricg$ and metric $\metricg$ on $Q$,
\begin{equation}\label{eq:bnablatildeg} (\dualconnection \cometricg)_X(\onef_1,\onef_2) := L_X \cometricg(\onef_1,\onef_2) - \cometricg(\dualconnection_X \onef_1, \onef_2) - \cometricg( \onef_1, \dualconnection_X \onef_2)
,
\end{equation}
describing the extent to which $\cometricg$ differs from the dual metric $\dualmetricg$.
This tensor will be the geometric source of ``bias" in the solutions of our control problem associated with $(Q,g,\tilde g)$, see Section~\ref{sec:biased}. In
Appendix~\ref{app:coordinates} we review a coordinate formula for the tensor $(\nabla^*\tilde g)$.

\subsection{The torque-cost cometric}\label{sub:tg}
In this subsection, we discuss in detail a particular type of cometrics $\tilde g$ which we call \emph{torque-cost cometrics} $\tg$. As explained in the introduction, these cometrics appear naturally in optimal control problems related to robotics and provide the main motivation for our study. The aim of this subsection is to detail practical characterizations of this cometric $\tg$ in relation to the mechanichal system's definition.

\smallskip

A typical robot can be modeled as a set of links connected by joints, with a set of \emph{actuators} (e.g., electric motors) attached to the mechanism and providing force or torque at specific points~\cite{LynchPark}. A robot is considered \emph{fully actuated} if the number of (independent) actuators is equal to the dimension of the configuration manifold, such that any desired force cotangent to the manifold can be produced by the actuators. A reasonable model for the cost of producing these forces is the sum of squared motor torques, which captures the property that ``overhead" costs in the system (e.g. heating of electrical components) scale super-linearly with the torque being produced.\footnote{More complete accountings of the cost of applying torques are available, accounting for properties such as the efficiency with which the system can regeneratively store energy while applying a joint torque in a direction opposite to that joint's motion~\cite{Abate:2016aa,TitusSpenny}, but for the purposes of this paper, we will restrict our attention to the torque-squared measure of effort.} We now outline how such torque-cost cometric $\tg$ is defined.

First, recall that the inertial metric $\metricg$ for such a robot is typically inherited from the Euclidean metric on its constituent particles through $Q\hookrightarrow Q_\amb$, as discussed in \S\ref{sec:kinematics}. This $g$ can be used to calculate the covariant acceleration $a$ of the robot as it executes a trajectory as well as to map this acceleration to the system force required to produce it. The dual of this metric, $\dualmetricg$, can be used to calculate a norm of the force, but the physical meaning of this norm is not particularly useful for the robot system: It describes the \emph{ambient-space magnitude} of the non-constraint force acting on the system.%

The ambient-space and actuator force magnitudes are not consistently proportional to each other, because the ambient-to-actuator force mapping for a rotary joint depends on the instantaneous moment-arm between the force and the joint. This means that to find the actuator effort corresponding to a trajectory we must \emph{incorporate this additional information about how the actuators are physically arranged in our system}. Thus, it is more useful to define our control problems %
in terms of a cometric $\tg$ that uses the induced actuator torques
\begin{equation}
\text{forces } F\mapsto (\tau(F)_1,\dots,\tau(F)_d) \text{ torques implemented by the actuators.}
\end{equation}
Typically, the underlying cost is a weighted sum of squared individual actuator torques, so that we obtain
\begin{equation} \label{eq:torquemetricdef}
\tg(F,F) = \sum_{i} k_{i}\bigl(\tau_{i}(F)\bigr)^{2},
\end{equation}
where the weighting terms $k_{i}$ encode the relative costs of generating torques with the actuators.

The following example illustrates how the actuator's configuration information is non-trivially encoded in $\tg\neq g^*$.

\begin{example}
Consider a point mass attached to a massless two-link arm as illustrated in Fig.~\ref{fig:torquesetup}. Taken as a point mass in $\euclid^{2}$, its mass matrix and dual mass matrix are
\begin{figure}[t]
\centering
\includegraphics[width=.7\textwidth]{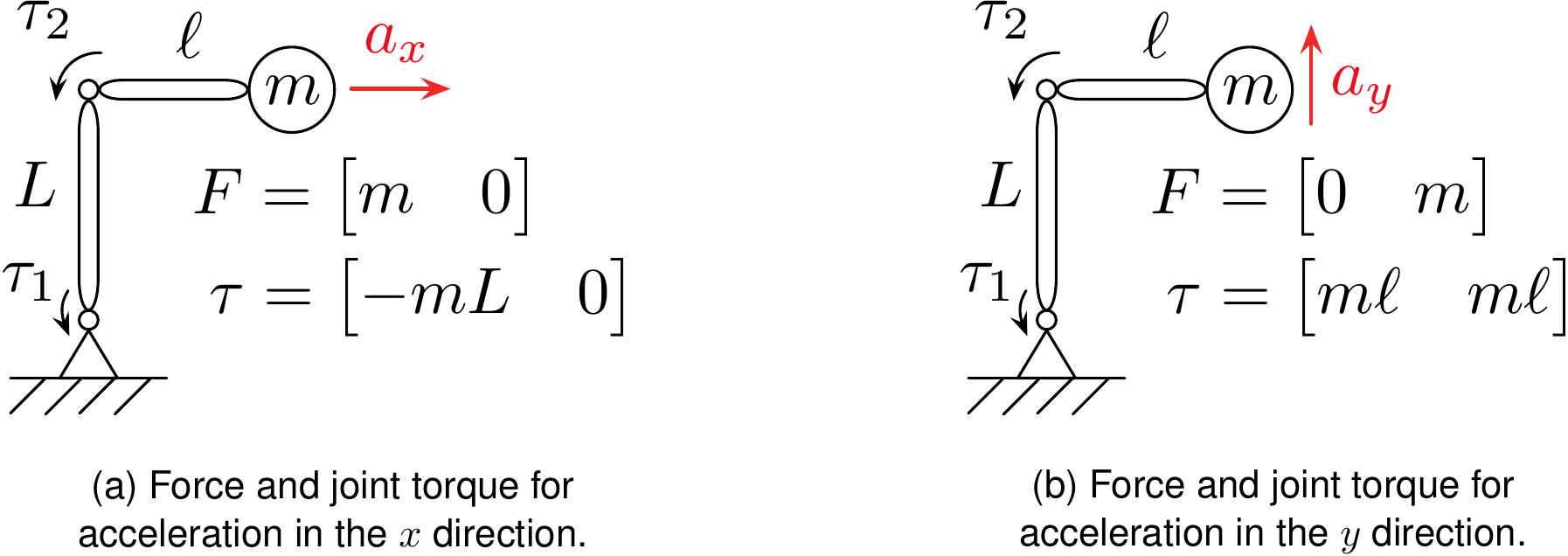}
\caption{The relationship between the forces felt by the system masses and those exerted by the system actuators depends on the relative position of the actuators to the masses and the directions in which the masses are (covariantly) accelerating.}
\label{fig:torquesetup}
\end{figure}
\begin{equation}
\metric = \begin{bmatrix} m \end{bmatrix} \qquad \mathrm{and} \qquad \dualmetric = \inv{\metric} = \begin{bmatrix} \frac{1}{m} \end{bmatrix}.
\end{equation}
The forces associated with unit acceleration in the $x$ and $y$ directions are
\begin{equation}
F_{x} = \metricg(\transpose{\begin{bmatrix} 1 & 0 \end{bmatrix}},\bullet) = \begin{bmatrix} m & 0 \end{bmatrix} \qquad \mathrm{and} \qquad F_{y} = \metricg(\transpose{\begin{bmatrix} 0 & 1 \end{bmatrix}},\bullet) = \begin{bmatrix} 0 & m \end{bmatrix},
\end{equation}
and the norms of the forces under the $\dualmetricg$ dual metric are equal to each other,
\begin{equation}
\norm{F_{x}}^{2}_{\dualmetricg} = \begin{bmatrix} m & 0 \end{bmatrix} \begin{bmatrix} \frac{1}{m} \end{bmatrix} \begin{bmatrix} m \\ 0 \end{bmatrix} = m \qquad \mathrm{and} \qquad \norm{F_{y}}^{2}_{\dualmetricg} = \begin{bmatrix} 0 & m \end{bmatrix} \begin{bmatrix} \frac{1}{m} \end{bmatrix} \begin{bmatrix} 0 \\ m \end{bmatrix} = m.
\end{equation}

If we instead take the system as a linkage, the Jacobian from the joint angles to the point mass motion in the illustrated configuration is
\begin{equation}
\jac = \begin{bmatrix} -L & 0 \\ \pmi \ell & \ell \end{bmatrix},
\end{equation}
such that the joint torques (calculated as the pullbacks of the $F_{x}$ and $F_{y}$ forces through the Jacobian) are
\begin{equation}
\tau_{x} = \begin{bmatrix} m & 0 \end{bmatrix} \begin{bmatrix} -L & 0 \\ \pmi \ell & \ell \end{bmatrix} =
\begin{bmatrix}
-mL & 0
\end{bmatrix}
\qquad \mathrm{and} \qquad
\tau_{y} =
\begin{bmatrix} 0 & m \end{bmatrix} \begin{bmatrix} -L & 0 \\ \pmi \ell & \ell \end{bmatrix}
=
\begin{bmatrix}
m\ell & m\ell
\end{bmatrix}
\end{equation}
These joint torques are clearly distinct from each other, depending on the location of the joints relative to the mass and the direction in which it is accelerating.

If we construct the system's mass matrix and its dual in the joint coordinates,
\begin{equation}
\metric = \transpose{\jac} m \jac = m\begin{bmatrix} (L^{2}+\ell^{2}) & \ell^{2} \\ \ell^{2} & \ell^{2} \end{bmatrix} \qquad \mathrm{and} \qquad \dualmetric = \inv{\metric} = \frac{1}{mL^{2}} \begin{bmatrix} \pmi 1 & -1 \\ -1 & (L^{2}+\ell^{2})/\ell^{2} \end{bmatrix},
\end{equation}
however, we find that the norms of both torques under $\dualmetricg$ are the same and are equal to the norm of the force,
\begin{align}
\norm{\tau_{x}}^{2}_{\dualmetricg} &=
\begin{bmatrix}
-mL & 0
\end{bmatrix}
\frac{1}{mL^{2}} \begin{bmatrix} \pmi 1 & -1 \\ -1 & (L^{2}+\ell^{2})/\ell^{2} \end{bmatrix}
\begin{bmatrix}
-mL \\ 0
\end{bmatrix}
= m
\\
\norm{\tau_{y}}^{2}_{\dualmetricg} &=
\begin{bmatrix}
m \ell & m \ell
\end{bmatrix}
\frac{1}{mL^{2}} \begin{bmatrix} \pmi 1 & -1 \\ -1 & (L^{2}+\ell^{2})/\ell^{2} \end{bmatrix}
\begin{bmatrix}
m\ell \\ m\ell
\end{bmatrix}
= m,
\end{align}
i.e., the dual metric $\dualmetricg$ uses information about the mechanism to \emph{remove} the specific geometric location of the actuators from the norm calculation, replacing the torque \emph{exerted} by the motors with the force \emph{felt} by the masses.

If we take the cost of supplying these torques as being individually intrinsic to the motors at the joints, therefore, we must (as detailed below) introduce a cometric $\cometricg$ distinct from the dual metric $\dualmetricg$ that captures these actuator efforts.

\end{example}

\subsubsection{Constructing the torque cometric}
Theoretically, we can consider $\tilde g$ as a given piece of information.
In practice, one way to construct a system's cometric $\tilde g$ is to select an atlas of preferred coordinate charts for $Q$ over which the underlying costs $\tilde g(F,F)$ are easily given, and then compute how its associated tensor $\cometric$ appears in those coordinates.
We now discuss this procedure for the specific case of the torque cometrics $\tg$ and in particular cases of interest.

In the case of the configuration space of a serial chain, it typically admits angular coordinates
\begin{equation}
\configspace \simeq \sphere^1\times \dots \sphere^1\times I \times \dots \times I,
\end{equation}
in which the $\sphere^1$ subspaces correspond to the unrestricted rotary joints, and the intervals $I\subset \R$ correspond to restricted rotary or sliding joints. Then, configurations can be parameterized (up to full rotations of the unrestricted joints, which can be handled by considering the universal cover of the corresponding $n$-torus) with a single chart, in which the chart parameters encode the joint angles or extensions. We shall refer to these as \emph{joint coordinates}. We then note that in the cotangent basis induced by this chart, the coefficients in the joint force vector $F$ are the individual joint torques,
\begin{equation}
F = \begin{bmatrix} \tau_{1} & \ldots & \tau_{n} \end{bmatrix}.
\end{equation}
We can then construct the torque cometric by defining its tensor in the chart-induced basis as
\begin{equation} \label{eq:weightedTM}
\cometric_{\tau} = \diag{(k_{1},\ldots,k_{n})},
\end{equation}
such that the force-norm under this metric meets our original definition,
\begin{equation}
\tg(F,F)   = \sum_{i} k_{i}\tau_{i}(F)^{2}.
\end{equation}

\begin{example}
A particular illustrative case arises when the costs for generating torques in the joints are all equal. Then, the weights can be taken as having uniform unit value, $k_i=1$, and the tensor for the torque cometric reduces to an identity matrix,
\begin{equation} \label{eq:IDTM}
\cometric_{\tau} = \matrixid^{n\times n}.
\end{equation}
In this case, the acceleration metric $h$ induced by $g$ and $\tilde g$ as in eq. \eqref{eq:defhfromtg} %
has a metric tensor given by
\begin{equation} \label{eq:combinedmetric}
\inmetric_{\tau} = \metric \cometric_{\tau} \metric = \metric\, (\matrixid^{\qdim\times \qdim})\, \metric = \metric^{2},
\end{equation}
the square of the inertial metric tensor as expressed in the joint-angle coordinates. The torque effort corresponding to a given covariant acceleration is thus
\begin{equation}
E = \inmetricg_{\tau}(a,\bullet) = \inmetric_{\tau} a = \metric^{2} a.
\end{equation}
\end{example}

\smallskip

\noindent {\bf Coordinate transforms of the torque cometric tensor.}
We now discuss briefly the glueing of the special-coordinate characterization of $\tg$ to yield a global tensor on $Q$.
The tensor for the torque cometric can be pulled back into coordinate charts other than those of the joints as
\begin{equation} \label{eq:cometrictransformation}
\cometric'_{\tau} = \transpose{J_{\text{t}}} \cometric_{\tau}  J_{\text{t}},
\end{equation}
where $J_{\text{t}}$ is the Jacobian of change from joint coordinates to the new ones. The tensor $M$ for the metric $g$ on $Q$ pulls back to the new coordinate chart as
\begin{equation}
M' = \invtranspose{J_{\text{t}}} M \inv{ J_{\text{t}}},
\end{equation}
so that the induced metric $h$ defined by eq.~\eqref{eq:defhfromtg} has the following tensor in the new coordinate bases
\begin{equation} \label{eq:antidualnewcoordinate}
\inmetric_{\tau}' = M' \cometric'_{\tau} M' = (\invtranspose{J_{\text{t}}} M \inv{ J_{\text{t}}})( \transpose{J_{\text{t}}} \cometric_{\tau}  J_{\text{t}})( \invtranspose{J_{\text{t}}} M \inv{ J_{\text{t}}}).
\end{equation}

\smallskip

\noindent{\bf Cometric on the ambient space.}
If the manifold metric $\metricg$ is inherited from the metric on an ambient space cut by holonomic constraints through $Q\hookrightarrow Q_\amb$ and $J$ is the corresponding Jacobian~\eqref{eq:defJac}, then $\cometric_\amb = J \cometric \transpose{J}$ defines an induced family of degenerate bilinear forms on the cotangent spaces $T^*\configspace_{\amb}|_Q$. These forms can be described using adapted coordinates to the embedding $Q\hookrightarrow Q_\amb$ as a pullback of $\cometric$ into the ambient space so that the kernel of $J \cometric \transpose{J}$ is conormal to the $Q$ and that any pullback of a force $F$ on $Q$ into $Q_\amb$ lies entirely in the orthogonal complement to this kernel.

Additionally, if there does happen to be a physical meaning to the ambient space such that forces $F$ originate from corresponding forces $F_\amb$ on the ambient space (e.g., because the system is a mechanical structure subject to external forces), the norm
\begin{equation}
\transpose{F_{\amb}} \cometric_\amb F_{\amb}
\end{equation}
projects out any component of $F_{\amb}$ that is absorbed by the constraint forces and then computes the length of the corresponding force vector in the priviledged coordinate system.

\begin{remark} {\textsc{Matrix eversion}}
There is an interesting relationship between the coordinate representation of an ambient-Euclidean metric and the ambient representation of a torque-based cometric with uniform actuator costs,
\begin{equation}
M_{\coord} = \transpose{J} J \qquad \text{and} \qquad \cometric_{\amb} = J\transpose{J},
\end{equation}
To the best of our knowledge, matrix relationships of this form have not been formally named in the literature. We suggest that it should be described as ``$\cometric_{\amb}$ is the \emph{everse} of $\metric$ with respect to the choice of orthonormal basis", where we have selected the term ``everse" to parallel the sense of ``inverse" while capturing the idea that $\metric$ has been ``everted" (turned inside out) to produce $\cometric_{\amb}$.
\end{remark}

\subsection{Examples of torque-cost cometrics}

As a concrete demonstration of how the torque cometric works, we have prepared a sequence of three examples. The first example builds on our sphere example, the second generalizes the sphere example to more completely show the geometric principles at play, and the final set of examples illustrates the effects of using different actuator attachments on the same mechanism.

\begin{example}
\textsc{(Torque-cost cometric for a yaw-pitch system on the sphere)}
Consider a point mass attached by a massless link to a yaw-pitch mechanism as in Figure~\ref{fig:yaw-pitch}. The link and joints constrain the mass to move over the surface of a sphere, so the kinematic data of this system $(\configspace,\metricg_{\sphere})$---e.g., the kinetic energy and covariant acceleration associated with its motion---are the same as those in our previous examples (see Example~\ref{ex:massmatrixsphere}), and the joint angles for yaw and pitch exactly correspond to the longitude-latitude coordinates $q=(\lon,\lat)$.

If we were to exert forces on this mechanism directly on the point mass (e.g., by fixing orthogonal ``rocket engines" to the mass, pointing in directions perpendicular to the link), then it would also be appropriate to measure the cost of such forces using the mass dual metric $\dualmetricg_{\sphere}$.
\begin{figure}
\centering
\includegraphics[width=.3\textwidth]{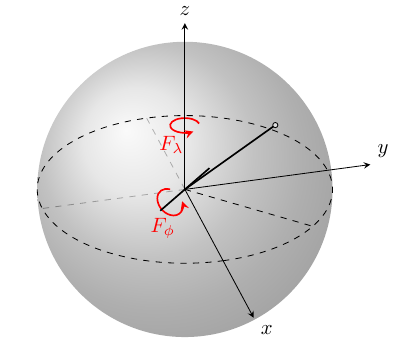}%
\includegraphics[width=.3\textwidth]{cometric5.pdf}%
\includegraphics[width=.3\textwidth]{cometric5.pdf}%
\caption{Given a yaw-pitch mechanism, the two most ``natural" ways of specifying the force are as covectors tangent to the sphere defined by the end of the mechanism, or as the torques around the two joints.}
\label{fig:yaw-pitch}
\end{figure}
If we instead apply forces to the system as torques around the yaw and pitch axes, then it is more appropriate to consider a torque-cost cometric $\tg$ which takes standard normal form in the yaw-pitch coordinates, as illustrated in the middle-right panel of Fig.~\ref{fig:spheremetricandtorquecometric}. %

In those coordinates, following eq.~\eqref{eq:combinedmetric}, the force norm $\norm{a}_{(\metricg,\cometricg)}$ is defined by the matrix
\begin{equation}
\inmetric_{\sphere\tau} = \metric_{\sphere}^{2}(\config) = \begin{bmatrix}
\cos^{4}{\lat} & 0 \\ 0 & 1
\end{bmatrix},
\end{equation}
as illustrated in the lower-right panel of Fig.~\ref{fig:spheremetricandtorquecometric}.
\begin{figure}[p]
\centering
\includegraphics[width=\textwidth]{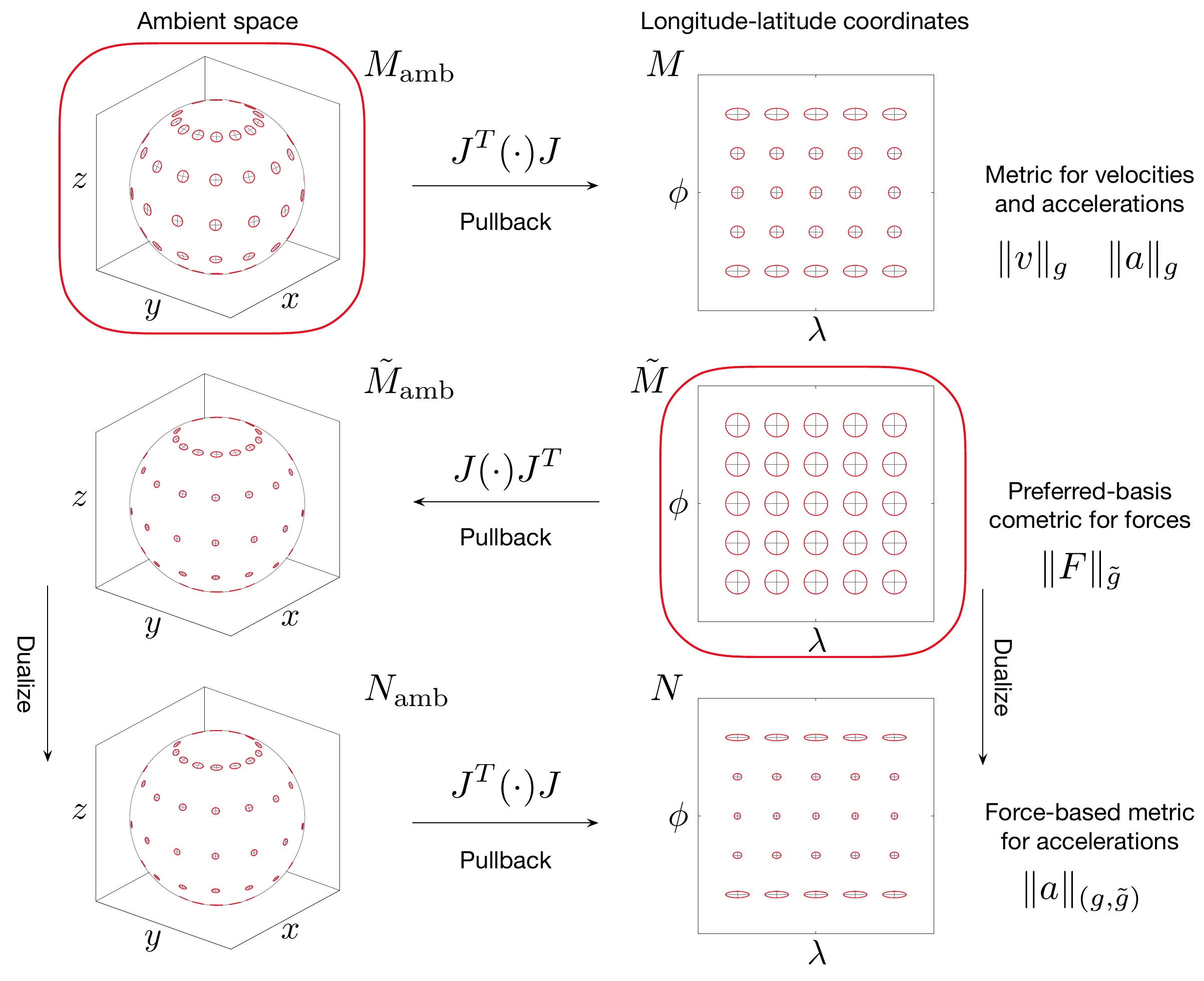}
\caption{Metric, preferred-basis cometric, and force-based acceleration metric for the sphere with $\lon,\lat$ parameterization. As before, the dualization operations are defined by the standard Euclidean ``identity" dualizations on the ambient space, and can be constructed on the coordinate chart by composition with the relevant pullback operations.}
\label{fig:spheremetricandtorquecometric}
\end{figure}
With respect to the $\metricg$-induced geometry at the left of Fig.~\ref{fig:spheremetricandtorquecometric}, the cometric $\tg$ indicatrices are stretched on the sphere---in configurations where the pitch is large, the yaw motor can exert more ambient force per unit torque, because the moment arm from the axis to the point is small. Because the longitude-lattitude coordinates define natural orthogonal directions for both $\metricg$ and $\tg$, the indicatrices for the combined metric $\norm{a}_{(\metricg,\tg)}$ are also streched similarly to those of $\norm{a}_{(\metricg)}$, but with quadratic intensity due to the fact that the underlying matrix is $\metric_{\sphere}^2(q)$.

From the ``effort" perspective, the mapping $E = \metric_{\sphere}^2(q)\, a$ captures the property that longitudinal accelerations near the pole are doubly-inexpensive to produce, because they are small covariant accelerations and because the moment arm from the yaw motor to the mass is small.

\end{example}

\begin{example}
\textsc{(Torque-cost cometric for a toroidal yaw-pitch system)}

If we place a link of length $\ell$ relative to the original link between the yaw and pitch joints of the mechanism in Fig.~\ref{fig:yaw-pitch}, the configuration space of the attached mass becomes a torus immersed in $\euclid^{3}$ as in Fig.~\ref{fig:twolink} (and we can consider the sphere from Fig.~\ref{fig:yaw-pitch} to be a degenerate special case of this immersion). The metric $\metricg_{\torus}$ produced by restricting the ambient Euclidean metric to the immersed torus has a metric tensor
\begin{equation}
M_{\torus} = \begin{bmatrix} (\ell+\cos{\lat})^{2} & 0 \\ 0 & 1 \end{bmatrix},
\end{equation}
and the torque-based cometric tensor remains an identity metric in the yaw-pitch coordinates,
\begin{equation}
\cometric_{\tau} = \matrixid^{2\times2}.
\end{equation}
The induced metric tensor in these coordinates thus yields an anti-dual to the torque cometric
\begin{equation}
\inmetric_{\torus\tau} = M_{\torus}^{2} = \begin{bmatrix} (\ell+\cos{\lat})^{4} & 0 \\ 0 & 1 \end{bmatrix}.
\end{equation}
Shifting our attention from the spherical mechanism of the previous example to the present more general two-link system highlights several important aspects of the torque-based cometric:

\begin{figure}
\centering
\includegraphics[height=2in]{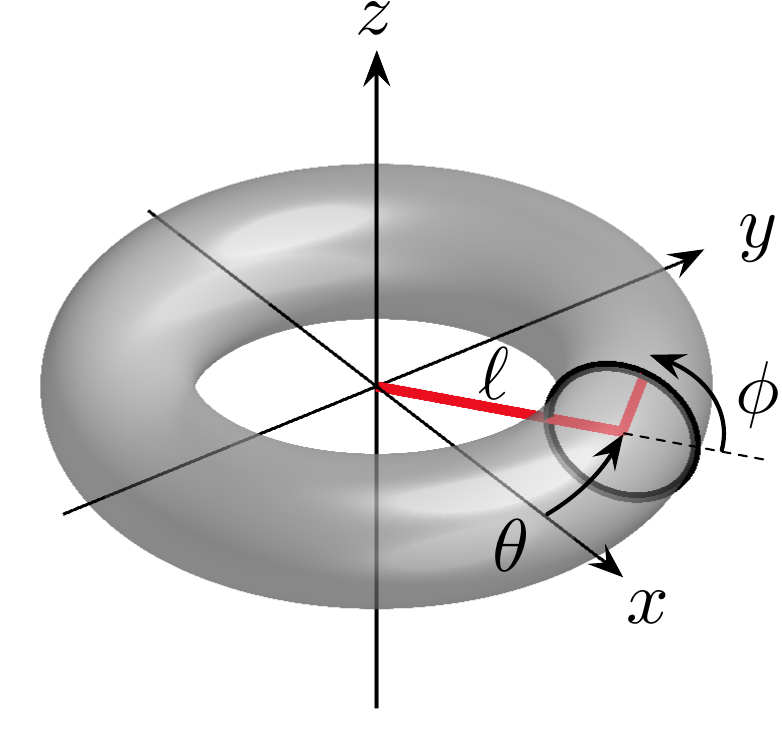}
\caption{Inserting a link of length $\ell$ between the yaw and pitch joints of the spherical mechanism from Fig.~\ref{fig:yaw-pitch} makes the configuration space of the end of the mechanism toroidal.}
\label{fig:toroidal}
\end{figure}

First, the configuration space of a mechanism conmposed of two unrestricted rotational joints is fundamentally toroidal.

Second, our construction of $\inmetric_{\torus\tau}$ represents here an interaction between the two ``natural" geometries associated with a torus: The metric $\metricg_{\torus}$ and its tensor $M_{\torus}$ are generated from the geometry of a torus immersed in Euclidean space (for which the length of a major-axis orbit depends on the minor-axis position at which the orbit is executed). In contrast, the torque-based cometric $\tg$ and its tensor $\cometric_{\tau}$ are derived from the geometry of a ``flat torus" ($\euclid^{2}$ with periodicity), because the cost of producing torque in a motor is generally independent of the present angle of the joints.

Finally, introducing the flat-torus geometry via the cometric does not ``blunt" the effects of the immersed-torus metric geometry. Rather, by making forces at large $\lat$ \emph{less} expensive than they would be under the ambient-force dual metric $\dualmetricg_{\torus}$, the torque-based cometric \emph{compounds with} the property of $\metricg_{\torus}$ that motions at large $\lat$ require less covariant acceleration.
\end{example}

\begin{example}
\textsc{(Torque-cost cometric for a planar mechanism actuated in two ways)}
Our final example systems, presented here together, illustrate the way in which different actuator arrangements interact with a mechanism to produce different $\cometric_{\tau}^{\antidual}$ matrices. The structure of both systems is a two-link planar chain with a point mass at its end, as illustrated in Fig.~\ref{fig:twolink}.

One system is actuated ``serially", with an actuator associated with each joint angle, so that the orientations of the links are $\theta_{1} = \alpha_{1}$ and $\theta_{2} = \alpha_{1}+\alpha_{2}$; the second is actuated ``in parallel", with an actuator associated with each link orientation, such that  $\theta_{1} = \beta_{1}$ and $\theta_{2} = \beta_{2}$.\footnote{This parallel actuation can be implemented physically by placing the $\beta_{2}$ motor axis coaxially with the $\beta_{1}$ motor axis, and then running a belt or other transmission element along the first link from the $\beta_{2}$ axis to a sprocket coaxial with the joint and fixed to the second link. The parallel actuation here does not produce any indeterminacy in mapping forces acting on the system to actuator forces (and so should not be confused with the closed-structural-loop construction of a \emph{parallel mechanism}, which can prevent there from being a deterministic mapping from system forces to joint forces).} The mapping from these actuator spaces onto the $xy$ position of the mass is illustrated in Fig.~\ref{fig:twolink} via sets of lines in which one actuator value is held fixed and the other is swept through its space of possible values.%

\begin{figure}
\centering
\includegraphics[width=\textwidth]{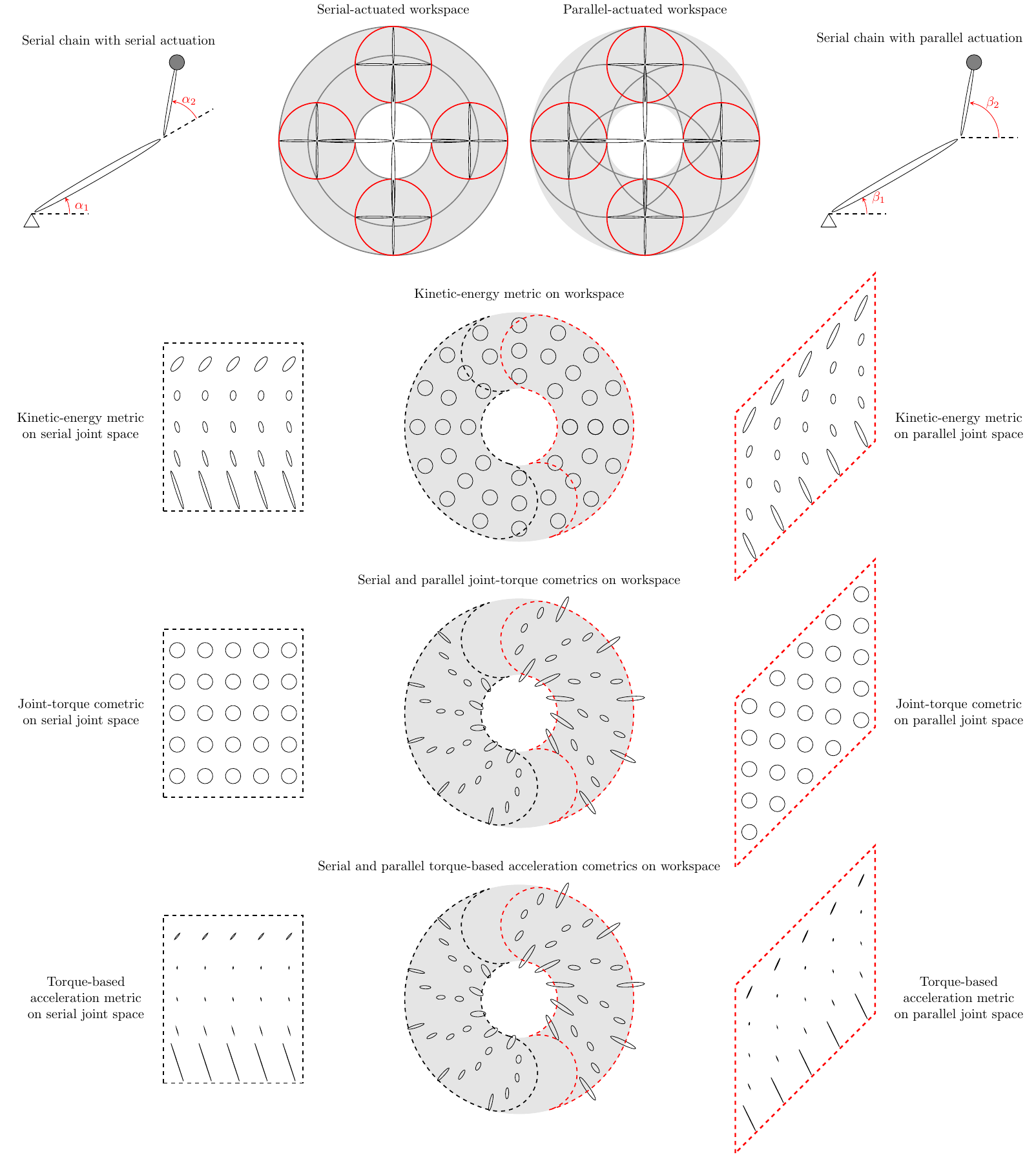}
\caption{Geometry of serial- and parallel-actuated two-link arms along with their metrics and cometrics.}
\label{fig:twolink}
\end{figure}

The kinetic-energy metric $\metricg$ is flat in the ambient $xy$ coordinates, and pulls back to the $\alpha$ and $\beta$ coordinates as illustrated in the second row of Fig.~\ref{fig:twolink}. The torque cometrics $\tg^{\alpha}$ and $\tg^{\beta}$ are flat in their respective coordinate spaces and pull back onto the ambient space with different values, as illustrated in the third row of Fig.~\ref{fig:twolink}.

As illustrated in the final row of Fig.~\ref{fig:twolink}, the flat nature of the torque cometrics means that the $\inmetric_{\tau}$ matrices of the system are the squares of the $\metric$ matrices on their respective spaces. Squaring the metric tensor does not commute with changes of coordinates,
\begin{equation}
(\transpose{\jac} M \jac)^{2} = \transpose{\jac} M \jac\transpose{\jac} M \jac \neq \transpose{\jac} M^{2} \jac,
\end{equation}
so these squared metric tensors represent truly different force-based norms of covariant acceleration. The fundamental difference between the two norms can also be seen directly in the $xy$ space, where the flat nature of $\metricg$ means that the two different pullbacks of $\cometric_{\tau}$ onto the ambient space also directly serve as the $\inmetric_{\tau}$ tensors in that space.
\end{example}

\clearpage
\section{Torque-optimal trajectories: biased geodesics and splines}\label{sec:biased}

In this section, we consider mechanical optimal control problems, with configuration space $Q$ and kinetic metric $g$, in which the system's forces are implemented by actuators and, thus, the optimization is relative to a general cometric $\cometricg$ abstracting the torque-cost cometric $\tg$ of Section~\ref{sub:tg}.

In the first subsection~\ref{sub:opcontrolprob}, we review the formulation of the optimal control problem associated with the data $(Q,g,\tilde g)$ and discuss a simple particular case in which the solution can be described by a 4th order ODE (see Example~\ref{ex:ELeq}). Already this example shows that, when the cometric is not the dual of the metric (i.e., $\tilde g\neq g^*$), the solution deviates from ordinary $g$-splines and $g$-geodesics.

Our main mathematical result in Section~\ref{sub:PMP} states that solutions of the general class of optimal control problems defined by $(Q,g,\tilde g)$ deviate from geodesics or splines for $(\configspace,\metricg)$ by intrinsic geometric effects involving the cometric $\cometricg$. The solution curves are described by an explicit Hamiltonian system of equations with $4dimQ$ unknowns or, equivalently, by a 4th order ODE for the trajectory $q(t)\in Q$ generalizing the Riemannian cubic $g$-spline equations \cite{CrLe,NoHePa,BCKS} (which we recover when $\tilde g=g^*$). We will call the general optimal trajectories {\bf biased splines and geodesics} associated with the input data $(\configspace,\metricg,\cometricg)$.

In Section~\ref{sub:descbiased}, we describe conceptually these biased curves and contrast them with the standard ones. Finally, in Section~\ref{sub:exsbiased}, we provide illustrative examples of these control problems and their solutions.

\subsection{Optimal control problems associated with $\mathbf{(\configspace,\metricg,\cometricg)}$}
\label{sub:opcontrolprob}

The control problems underlying the interests of this paper (e.g., moving a robot between two configurations) can be stated more abstractly as finding trajectories $q(t)\in Q$ connecting initial and final states in a given time at a minimum of control effort (see e.g. \cite{BulLew,FraMa,GaHoMeRaVi1}). Taking the control effort as an appropriate square of the system force, these control problems can be formally stated as optimization problems with
\begin{subequations}
\begin{align}
\text{State:} && x &: t \mapsto \begin{bmatrix} \config(t) \\ \configdot(t) \end{bmatrix} \in T\configspace \label{eq:optstate}\\
\text{Control:} && u&= F = \metricg\bigl(\nabla^{\metricg}_{\dot{\config}}\dot{\config},\bullet\bigr)
\label{eq:optcontrol} \\
\text{Cost functional:} && I &  = \int_{t_0}^{t_f} \cometricg(F,F) \ dt \label{eq:control}
\end{align}
in which the goal is to solve for a state trajectory $x(t)\in T\configspace$ that minimizes the cost functional over the period of the motion while satisfying the boundary-condition constraints
\begin{equation}
(q(t_0),\configdot(t_0)) = (q_{0},\configdot_{0}) \qquad \text{and} \qquad (q(t_f),\configdot(t_f)) = (q_{f},\configdot_{f}) \label{eq:optconstraints}
\end{equation}
\end{subequations}
on the endpoint positions and velocities of the trajectory. Note that the force is defined through the covariant acceleration $a(\config,\configdot,\ddot{\config}) = \nabla^{\metricg}_{\configdot}\configdot$ only using the mass metric $\metricg$, whereas the cost functional uses an arbitrary cometric $\cometricg$.

The optimal control problem~\eqref{eq:control} is equivalent to one in which we consider the acceleration $a$ as being the control parameter instead of the force, $F= g(a,\bullet)$. This problem will be easier to handle at the theoretical level as it is closer to previously studied systems (in particular to those studied in~\cite{BCKS}) and we now proceed to describe it.

The acceleration-based control problem is defined by the same states $(q,\configdot)\in TQ$ defined by a curve $q(t)\in Q$ together with
\begin{eqnarray}
\text{Control:} \ && u'= a = \nabla^g_{\configdot} \configdot  \nonumber \\
\text{Cost functional:} \ && I' = \int_{t_0}^{t_f} h (a, a) \ dt = \int_{t_0}^{t_f} \tilde g (g(a,\bullet),g(a,\bullet)) \ dt  \label{eq:controlh}
\end{eqnarray}
where $h$ is defined from $\cometricg$ and $g$ as in eq.~\eqref{eq:defhfromtg}. The boundary conditions are the same as in~\eqref{eq:optconstraints}. It is clear that $q_{\min}(t)\in Q$ is a solution of the optimal control problem~\eqref{eq:control} iff it also solves~\eqref{eq:controlh}. In this case, the optimal controls, $F_{\min}$ and $a_{\min}$, for the underlying solution $q_{\min}(t)$ are related via $F_{\min}=g(a_{\min},\bullet)$.

\begin{remark} \emph{(Time intervals)}
Time-symmetry allows us to introduce two simplifcations to the optimal control problem: First, we can start the ``clock" for the system at the beginning of the trajectory, such that $t_{0} = 0$. Second, because rescaling time introduces a constant scaling factor into the acceleration (and thus actuator-force) calculations, we can restrict our attention to motions of unit time-period, such that the interval $[t_{0}, t_{f}]$ becomes $[0,1]$ (noting that any velocities in the specified boundary conditions should be rescaled to match the new timescale).
\end{remark}

\noindent {\bf Discussing approaches to the solution.} In solving optimal control problems of this form, there are two fundamental strategies: {\bf shooting} approaches and {\bf variational} approaches. (For general variational methods on manifolds, see e.g.~\cite{GaHoMeRaVi1,GaHoMeRaVi2}.) Both are fundamentally based on Pontryagin's maximum principle (PMP, see~\cite{BuHoMe}), which describes the necessary conditions that an optimal trajectory must satisfy. These conditions correspond to admitting a lift to Hamiltonian system defined on $T^*(T\configspace)$ or, equivalently, defining  a fourth-order differential equation for the underlying system configuration $q(t)\in Q$.

The standard shooting approach based on PMP will be described in Section \S\ref{sub:PMP} below (see, in particular, Remark~\ref{rmk:PMPshooting} for a discussion of the shooting aspect). This will lead us to a complete theoretical characterization of the solutions in terms of differential equations.

Variational approaches can be also useful in practice. In these approaches, the optimal trajectories are found by considering \emph{search space} which is typically a finite dimensional family of curves that satisfy the boundary conditions~\eqref{eq:optconstraints} (e.g., the locations of a set of waypoints distributed along the cuve). Starting with an arbitrary curve in the search space, this curve is deformed until it satisfies the extended Hamiltonian differntial equation at its interior points.

The relative efficiencies of the shooting approach and the variational approach depend on how expensive it is to repeatedly evolve the Hamiltonian differential equation as compared to finding the gradient of Hamiltonian-satisfaction as a connecting trajectory is deformed. Hybrid approaches are also available; for example, ``multiple shooting methods" evolve the extended Hamiltonian equation from several starting points between the initial and final states, and then adjust the positions and initial conditions of each of these points so that the evolved trajectories form a unified whole.

\medskip

We finish this subsection with a simple formal example that illustrates the role of $\tilde g$ in the variational problem.
\begin{example}\label{ex:ELeq} \emph{(How $\tilde g$ enters the Euler-Lagrange equations)} Consider a simple (purely mathematical) situation in which $Q=\R^n$ and $g$ is the flat Euclidean metric. In this case, we can naturally identify $(\R^n)^*\simeq \R^n$ under which $F\simeq a=\ddot q$, the second ordinary time-derivative. The cost functional~\eqref{eq:controlh} yields
\begin{equation}
I' = \int_{t_0}^{t_f} h(\ddot q, \ddot q) \ dt
\end{equation}
where $h\equiv h_{ij}(q)$ is an arbitrary metric induced by an underlying arbitrary cometric $\tilde g$ as in~\eqref{eq:defhfromtg} (which here, because of the flat metric, simply implies lowering the indices by means of the Euclidean $g$). The boundary conditions on $(q,\configdot)\in \R^n \times \R^n$ are the same as in~\eqref{eq:optconstraints}. To get the corresponding Euler-Lagrange equations for this variational problem, we consider as usual a variation $\gamma\equiv q(t,\epsilon), \ \epsilon\in (-\delta,\delta)$, respecting the boundary conditions, and sample for critical points of $I'$,
\begin{equation}
\frac{\partial I'[\gamma]}{\partial \epsilon}|_{\epsilon =0} = 0.
\end{equation}
A straight-forward computation, involving twice integrating by parts, leads to the following 4th order ODE for the critical curve $q(t)\in Q$,
\begin{equation}\label{eq:ELeqsimple}
\frac{d^2}{dt^2}[h_{ij}\ddot{q}^j]=-\frac{1}{2}\partial_ih_{jk} \ddot{q}^j \ddot{q}^k, \ \forall i=1,\dots,n.
\end{equation}
Notice that when $h_{ij}(q)=h_{ij}$ is constant, we obtain the Euclidean cubic spline equation: $d^4 q/dt^4 = 0$. For general $h_{ij}(q)$, we see that its derivatives introduce a ``bias" with respect to the Euclidean case. In the following subsection we shall interpret the extra terms appearing in an intrinsic geometric way, even for general $(Q,g,\tilde g)$.
\end{example}

\medskip

\subsection{General solutions from Pontryagin's Maximum Principle}\label{sub:PMP}

In this subsection, we provide a system of differential equations which generalize to arbitrary triples $(Q,g,\tilde g)$ both the ordinary cubic spline equations associated with a $(Q,g)$ and the ones obtained above in the particular case of Example~\ref{ex:ELeq}. To handle the geometry in an invariant way, we resort to Pontryagin's Maximum Principle (PMP, see~\cite{BuHoMe}) and the techniques devolopped in~\cite{BCKS}.

We then focus on applying PMP to our optimal control problem~\eqref{eq:controlh}  defined by $(\configspace,\metricg,\cometricg)$. This principle asserts that the solution $(\gamma(t), \dot \gamma(t)) \in T\configspace$ is the projection under the natural  map $T^*(T\configspace) \to T\configspace$
of a solution of an extended Hamiltonian system $(T^*(T\configspace) , \omega_c, H_{\max})$. Here, $\omega_c$ is the canonical cotangent-bundle symplectic structure and the Hamiltonian function
\begin{equation}
H_{\max} = \underset{u'}{max}\{H_{u'}\}
\end{equation}
is obtained as the maximum over a family of hamiltonians parameterized by the control $u'$ (which in the case of~\eqref{eq:controlh} is the acceleration $u'=a$). Also notice that this extended system has twice as many variables with respect to our initial first order problem on $T\configspace$. The extra variables will be denoted $(\alpha,p)$ and are thought of as auxiliary `co-state' variables dual to the original $(q,v)$ in $T\configspace$, respectively (see more below).

\begin{remark}\label{rmk:PMPshooting}
Note that, in practice, one needs to find the correct initial conditions $(\alpha_0,p_0)$ for the additional co-state variables $(\alpha,p)$ in this extended Hamiltonian system so that the projected curve satisfies the desired final-time boundary conditions $(q_f,\configdot_f)$ of the original control problem: the optimal control problem is transformed into a shooting problem in a Hamiltonian system with extra auxiliary variables.
\end{remark}

In~\cite{BCKS}, such Hamiltonian systems $(T^*(T\configspace) , \omega_c, H_{\max})$ were analyzed using the connection $\nabla\equiv \nabla^g$ to split the variables as $T^*(T\configspace)\simeq T^*\configspace \times_\configspace T\configspace \times_\configspace T^*\configspace$ in order to find the optimal Hamiltonian and the underlying Hamiltonian equations explicitly. Nevertheless, only a restricted class\footnote{The functionals considered in \cite[Prop. 2]{BCKS} are of the form $\int^{t_f}_{t_0} c(g(a,a))\ dt$ with $c:\mathbb{R}\to \mathbb{R}$.} of cost functionals was considered in~\cite{BCKS} and we need to slightly generalize the arguments described in that reference to apply them to our control system defined by $(\configspace,\metricg,\cometricg)$. In this way, we get the following:

\begin{theorem}\label{thm:solutions}
A curve $\gamma(t)\in \configspace$ is a solution to the optimal control problem defined in~\eqref{eq:controlh} by $(\configspace,\metricg,\cometricg)$ iff there exist curves $\alpha(t)$ and $p(t)$ in $T^*_{\gamma(t)} \configspace$ such that $(q(t)=\gamma(t),v(t)=\dot \gamma(t),\alpha(t),p(t))$ is a solution of the following (Hamiltonian) system of differential equations:
\begin{eqnarray}
\configdot &=& v \nonumber \\
\nabla_{\configdot} v &=& \frac{1}{2} h^*(\alpha, \bullet ) \in T_{q(t)}Q \nonumber \\
\bnabla_{\configdot} \alpha &=& - p \nonumber \\
\bnabla_{\configdot} p &=& -\frac{1}{4} (\bnabla h^*)_{\bullet} (\alpha,\alpha) +\langle \alpha, R_g(\bullet,v)v \rangle \in T^*_{q(t)} Q, \label{eq:bsystemh}
\end{eqnarray}
where %
$h^*$ is the cometric defined as the dual of the $h$ given in eq.~\eqref{eq:defhfromtg}, $R_g\in \Omega^2(\configspace,End(T\configspace))$ is the Riemann curvature tensor of the Levi-Civita connection $\nabla\equiv \nabla^g$ associated with $g$ and $\bnabla h^*$ is the tensor~\eqref{eq:bnablatildeg} measuring the compatibility between the cometric $h^*$ and the connection $\bnabla$ on $T^*\configspace$ dual to $\nabla$.
\end{theorem}

\noindent {\bf Proof:} We need to verify that the equations coming from PMP, as briefly recalled above, coincide with the ones given in the theorem. As mentioned, we follow~\cite{BCKS} closely. First, the state equations $v=\configdot$ and $u'=\nabla^g_{\configdot}\configdot$ of our control problem are exactly those considered in that reference. The integrand in the cost functional~\eqref{eq:controlh} is given by the function $c(q,u')=h(u',u')|_q, q \in Q$ which is quadratic in the control. Follwing the general method of PMP, as explained in~\cite{BCKS} for the present type of systems, one can use the global ``splitting" diffeomorphism $T^*(TQ)\simeq T^*Q\times_Q TQ \times_Q T^*Q$ induced by the connection $\nabla\equiv \nabla^g$ and obtain the $u'$-family of Hamiltonians $H_{u',\nabla}:T^*Q\times_Q TQ \times_Q T^*Q \to \R$ given by
\begin{equation}
H_{u',\nabla} = -c(q,u') + \langle \alpha,u'\rangle + \langle p, v\rangle = -h(u',u')|_q +  \langle \alpha,u'\rangle + \langle p, v\rangle
\end{equation}
where $q\in Q$, $\alpha,p \in T^*_q Q$ and $v\in T_q Q$ as in the statement of the theorem. PMP then involves finding the maximum of this family $H_{u',\nabla}$ over the possible $u'$. In the present case, due to the quadratic behaviour, the solution is simple and given by $u'_{\max}=\frac{1}{2}h^*(\alpha,\bullet)$ where $h^*$ is the cometric dual to $h$. The maximized Hamiltonian is then given by
\begin{equation} \label{eq:Hmaxnabla}
H_{\max,\nabla} =\frac{1}{4} h^*(\alpha,\alpha) + \langle p, v\rangle.
\end{equation}
We thus need to compare the equations of motion for this Hamiltonian system $(T^*Q\times_Q TQ \times_Q T^*Q, \omega_\nabla, H_{\max,\nabla})$ with the ones provided in the theorem. To that end, first note that the symplectic structure $\omega_\nabla$ is the pullback of the canonical one on $T^*(TQ)$ under the mentioned splitting diffeomorphism and that it is described in detail in~\cite[Prop. 1]{BCKS}. Next, because the equations given in our theorem are globally defined, it is enough to verify that they coincide with the ones coming from PMP in a generic coordinate chart induced by coordinates on $Q$. Finally, following~\cite[Prop. 1(iii)]{BCKS}, such coordinate form of the PMP-induced equations of the Hamiltonian system
$(T^*Q\times_Q TQ \times_Q T^*Q, \omega_\nabla, H_{\max,\nabla})$ with the above explicit Hamiltonian function~\eqref{eq:Hmaxnabla} can be computed directly yielding:
\begin{eqnarray}
\configdot^i &=& v^i \nonumber \\
\dot v^a + \Gamma^a_{ib}\configdot^i v^b &=& \frac{1}{2} (h^*)^{ac}\alpha_c \nonumber \\
\dot \alpha_a - \Gamma_{ia}^b \configdot^i\alpha_b &=& -p_a \nonumber \\
\dot p_i - \Gamma_{ia}^b v^a p_b &=&-\frac{1}{4} \left( \partial_i(h^*)^{ab} \alpha_a \alpha_b + 2 \Gamma_{ia}^b (h^*)^{ac} \alpha_b\alpha_c \right) + R_{ija}^b \configdot^j v^a \alpha_b \label{eq:pmpsystcoord}
\end{eqnarray}
where all the indices $i,j,a,b,etc.$ run from $1,\dots,dimQ$, the Christoffel symbols of the connection are defined by $\nabla^g_{\partial_{q^i}} \partial_{q^a} = \Gamma_{ia}^b\partial_{q^b}$, and $R^b_{ija}$ denote the coefficients of the curvature tensor $R_g$ defined by $\nabla^g$. Appendix~\ref{app:coordinates} provides more details on these coordinate expressions and, in fact, using the coordinate version~\eqref{eq:coordhsystem} of~\eqref{eq:bsystemh} given in the Appendix, one readily verifies that the PMP-induced system~\eqref{eq:pmpsystcoord}  coincides exactly with the one given in the statement of the theorem. This finishes the proof. $\square$.

\medskip

The Theorem's system of equations~\eqref{eq:bsystemh} is expressed in terms of the relevant differential-geometric information underlying $(\configspace,\metricg,\cometricg)$.
An interesting particular case is that of a general $(Q,g)$ with cometric $\tilde g=g^*$ given by the dual metric. In that case, $h=g$ by~\eqref{eq:hequalsg} and, as mentioned before, the \ctensor\ vanishes, $\nabla^*g^*=0$, as a consequence of $\nabla^g g=0$. The equations can then be seen to reduce to the cubic Riemannian spline equations associated with $(Q,g)$, as described in~\cite{BCKS} (see also below). Thus, the term $\bnabla h^*$ can be interpreted as the geometric origin of the biasing away from the standard Riemaniann splines.

\begin{remark}\emph{(Some hidden ingredients)} As mentioned in the above proof, we remark that the symplectic structure $\omega_\nabla$ on $T^*\configspace \times_\configspace T\configspace \times_\configspace T^*\configspace$ which makes the system~\eqref{eq:bsystemh} Hamiltonian is explicitly described in~\cite{BCKS}, as well as the diffeomorphism taking it to the canonical one $\omega_c$ in $T^*(T\configspace)$. Also, we stress that the optimal control maximizing the underlying family of Hamiltonians $\{H_{u', \nabla}\}$ is explicitly given by $u'_{\max}=a_{\max}=\frac{1}{2} h^*(\alpha,\bullet) \in T_{\gamma(t)}\configspace$. Nevertheless, we shall not be using these ingredients further in this paper.
\end{remark}

\begin{remark} \emph{(Comparing to the variational equation for $(Q,g)$ Euclidean)}\label{rmk:eqshsol}
Notice that in the case $Q=\R^n$ and $g$ is the Euclidean metric, this 1st order (Hamiltonian) system~\eqref{eq:bsystemh} is equivalent to the 4th order ODE for $q(t)$ that we derived variationally in Example~\ref{ex:ELeq}. To see this, we use the first three equations to get that $\alpha=2h(\ddot q,\bullet)$ and $p = -\frac{d}{dt}(2h(\ddot q,\bullet))$ or, in standard coordinates, $\alpha_i = 2 h_{ij}\ddot q^j$ and $p_i= -\frac{d}{dt}(2h_{ij}\ddot{q}^j)$. We also recall that $R_g=0$ and that, since $h^*\equiv (h^*)^{ij}$ is given by the inverse matrix of $h\equiv h_{jk}$, we have the identity $\partial_i(h^*)^{jk} h_{kl} = - (h^*)^{jk} \partial_i h_{kl}$. Substituting the above into the last equation in~\eqref{eq:bsystemh}, we recover eq.~\eqref{eq:ELeqsimple} of Example~\ref{ex:ELeq}.
\end{remark}

It is straight-forward to verify that this Hamiltonian (1st order) system is equivalent to the 4th-order ODE for the curve $q(t)\in Q$ given by
\begin{equation}\label{eq:4orderh}
\bnabla_{\configdot}\bnabla_{\configdot}(h^\flat(\nabla_{\configdot}\configdot)) = -\frac{1}{2} (\bnabla h^*)_\bullet (h^\flat(\nabla_{\configdot}\configdot), h^\flat(\nabla_{\configdot}\configdot)) - h(\nabla_{\configdot}\configdot, R_g(\bullet,\configdot)\configdot),
\end{equation}
where we return momentarily to the standard notation $h^\flat:TQ\to T^*Q, v\mapsto h(v,\bullet)$ to avoid confusion with the other bullets appearing in the equation.
As observed in Remark~\ref{rmk:eqshsol}, this equation generalizes the one found variationally in Example~\ref{ex:ELeq} when $(Q=\R^n, g)$ was Euclidean  and, at the same time, organizes the extra $h$-terms (equiv. $\tilde g$-induced terms) in an intrinsic geometric way.
Moreover, when $\tilde g=g^*$, the above equation reduces to the 4th order cubic Riemannian spline equation in $(Q,g)$ (see \cite{CrLe,NoHePa} and~\cite[eq. (19)]{BCKS}, as well as eq.~\eqref{eq:riemannianspline} below). This follows from $\bnabla g^*=0$ and $\bnabla g^\flat = g^\flat \nabla$.

\subsection{Description of optimal trajectories} \label{sub:descbiased}
In this subsection, we describe conceptually the trajectories which solve the optimal control problem associated with $(Q,g,\tilde g)$. We begin with the known case in which $\tilde g=g^*$ and the trajetories are Riemanian splines or geodesics. We then use this case to describe by contrast what happens in the biased-splines and geodesics case (i.e. for general $\tilde g \neq g^*$).

\subsubsection{The dual metric case: Riemannian Splines and Geodesics}\label{subsec:splines}

If we take cometric $\tilde g=g^*$ to be the dual of $g$, then the induced metric $h$ defined by~\eqref{eq:defhfromtg} yields $h=g$. As mentioned before, the solutions to the control problem~\eqref{eq:controlh} (equiv.~\eqref{eq:control}) are given by cubic Riemannian splines on $(Q,g)$ (see~\cite{CrLe,NoHePa,BCKS}).
We shall provide in the sequel a conceptual analysis of these curves to be used later.

\smallskip

\noindent {\bf Review of cubic Euclidean splines.}
As mentioned in Example~\ref{ex:ELeq}, cubic splines in $\euclid^{n}$  are trajectories that can be described extrinsically as cubic polynomials in time,
\begin{equation}
q(t) = q_{0} + tq_{1} + t^{2}q_{2} + t^{3}q_{3}
\end{equation}
(giving them their name) or, equivalently, in terms of the intrinsic requirement of having zero snap,
\begin{equation}
s = \frac{d^{4}}{dt^{4}} x = \frac{d^{2}}{dt^{2}} \ddot{x} = 0.
\end{equation}
These curves are the lowest-order curves capable of connecting arbitrary position-and-velocity boundary conditions given by~\eqref{eq:optconstraints}.
As we saw in Example~\ref{ex:ELeq}, the cubic spline connecting a pair of position-and-velocity boundary conditions minimizes
$%
I'= \int_{t_{0}}^{t_{f}} |\ddot{q}|^2\, dt
$ %
over all such curves. Intuitively, this rests on the principle that the cost function favors trajectories that traverse short paths between the points (reducing the overall magnitude of the required acceleration),  and trajectories that spread their acceleration evenly over time (reducing the quadratic penalty for concentrating the acceleration at any one point). The linear acceleration profile of a cubic spline balances path shortness against curvature concentration.

\smallskip

\noindent {\bf Splines on manifolds.}
Acceleration-minimizing trajectories $\gamma(t)\in Q$ on manifolds with curvature or non-Cartesian parameterizations typically cannot be described extrinsically via cubic polynomials. It is, however, possible to generalize the intrinsic (snap-based) definition of a cubic spline onto the Riemannian case associated with $(Q,g)$. To describe this case, we first take the acceleration, jerk, and snap of the trajectory $\gamma(t)\in Q$ to be calculated covariantly along the trajectory, producing a system of equations
\begin{subequations}
\label{eq:covariantacceljerksnap}
\begin{align}
v &= \dot{\gamma} && \text{velocity}\\
a &= \nabla_{v}v && \text{covariant acceleration}\\
j &= \nabla_{v}a && \text{covariant jerk} \\
s &= \nabla_{v}j  && \text{covariant snap}
\end{align}
\end{subequations}
which can be integrated over time for known $(\gamma,v,a,j)$ initial conditions and a specified snap function.

The trajectory which minimizes $I'=\int_{t_0}^{t_f} g(a,a) \ dt$ is determined by  the condition (see e.g.~\cite[eq. (19)]{BCKS}) that the covariant snap of the acceleration is determined by the curvature $R_g$ as
\begin{equation} \label{eq:riemannianspline}
s = - R_g(a, \dot \gamma) \dot \gamma,
\end{equation}
for which the in-coordinate calculation formulas are given in Appendix~\ref{app:coordinates}.

Let us discuss conceptually the right-hand side of~\eqref{eq:riemannianspline}.
Intuitively, this curvature term describes how much ``free" velocity change the trajectory can capture by exploiting the curvature of the manifold:
In regions of the manifold with negative curvature,\footnote{Specifically, negative  curvature in the $(\dot \gamma, a_{\gamma})$ section. Note that the vectors defining the section appear in reverse order in the arguments of $R_{\metricg}$.} neighboring trajectories diverge from each other. Because this divergence amplifies the acceleration that separates the neighboring trajectories, a trajectory that ``invests" in early acceleration in a negative curvature region requires less total acceleration than a trajectory that keeps to a strictly linear acceleration profile, and so can have a lower total cost even though it pays a quadratic penalty on this initial investment.

Conversely, neighboring trajectories in regions of positive curvature converge; this convergence means that early accelerations are ``taxed" over the remainder of the trajectory and thus that it is optimal to delay acceleration to reduce the time over which this penalty is assessed (even though this delay incurs a penalty from the instantaneous quadratic cost).

The Riemannian curvature term itself measures the divergence or convergence of neighboring trajectories as the difference between the tangent vector an infinitesimal step along the geodesic the system is currently following and the tangent vector an infinitesimal step along the geodesic it is accelerating onto. The negative of this curvature term thus provides an answer to ``does the curvature of the manifold amplify or diminish the effect of acceleration applied at this point".

\smallskip

\noindent {\bf Force-based equation for ordinary $(Q,g)$ splines.}
For later comparison with the general $(Q,g,\tilde g)$ case, it is helpful to rewrite the above equations for the cubic Riemannian splines in terms of the equivalent force-based description with functional $I$ given by~\eqref{eq:control} with $\tilde g=g^*$. Using the relation $F=g(a,\bullet)$ and  denoting the first and second derivatives of the force as the ``yank"~\cite{LinYank} and ``tug"~\cite{JazarDynamics} experienced by the system, eqs.~\eqref{eq:covariantacceljerksnap} and~\eqref{eq:riemannianspline} become equivalent to
\begin{subequations} \label{eq:splineforce}
\begin{align}    v &= \dot{x} &&&&&& \text{velocity} \label{eq:splineforcevel}\\
a &= \nabla_{v}v && && &&  \text{covariant acceleration}\\
&& F &= \metricg(a,\bullet) && && \text{force is dual to acceleration} \label{eq:accelbulletdual}\\
&& && \yank &= \dualconnection_{v}F  &&  \text{covariant yank} \\
&& && \tug &= \dualconnection_{v}\yank =  \covprod{F}{-(R_g(\bullet, v) v)}, &&  \text{covariant tug} \label{eq:splineforcetug}
\end{align}
\end{subequations}
where $\dualconnection$ is the connection dual to $\nabla^g$ as before. Note that because we set $\tilde g=g^*$, we have $h=g$ (see~\eqref{eq:defhfromtg})  and $\nabla^* h^*=0$.

In the tug expression from~\eqref{eq:splineforcetug},  $\covprod{F}{-(R_g(\bullet, v) v)}$ measures the extent to which force applied to the system at the current time will ``pay dividends" (reduce the future force needed to achieve a given change in velocity) or ``be taxed" (such that the effect of any force applied to the system decays over time), and thus whether the system should ``invest" or ``withhold" force from the system. The construction of this term is somewhat more intricate than the simple dual of acceleration: because $R_{\metricg}(a,v)v$ returns a vector value and not a scalar, $R_{\metricg}(\bullet,v)v$ is a ``covector of vectors", which is reduced to a simple covector by the product operation of $F$ with each of these vectors.

The above system is then equivalent to the Hamiltonian system~\eqref{eq:bsystemh} under $\alpha=2 h(\nabla_{v}v,\bullet) = 2F$ and $p=-2 \yank$. This fact also recovers the Hamiltonian description obtained in~\cite{BCKS} of the 4th order $(Q,g)$-cubic spline equation~\eqref{eq:riemannianspline}.

\smallskip

\noindent {\bf Special case solutions.}
Under certain classes of start and end velocity constraints on the trajectory, cubic Riemannian splines have additional interesting properties:%
\begin{enumerate}
\item If the initial and final velocities $\configdot_{0}$ and $\configdot_{f}$ are left unspecified, then the solutions to the optimization problem are geodesic trajectories on $(\configspace,\metricg)$ between $q_{0}$ and $q_{f}$. To see this, observe that geodesic trajectories have no covariant acceleration, such that if there exists a geodesic passing through $q_{0}$ and $q_{f}$ (i.e., the start and end points lie in a geodesically-convex region of the configuration space), the geodesic is a minimizer of the squared covariant acceleration over the trajectory. If there is more than one geodesic connecting the points, all of these geodesics minimize the squared covariant acceleration to zero (but they can be discriminated amongst by noting that the time taken to traverse each geodesic at a given level of kinetic energy is proportional to its pathlength).
\item If the initial and final speeds are zero (i.e., the commanded trajectories are point-to-point motions with zero velocity at the ends), then the solutions to the optimization problem follow the same paths as the geodesics of $(\configspace,\metricg)$, but with position along the path described by a cubic function of time rather than the linear (constant-speed) progression of a geodesic. To see this, assuming again that $q_{0}$ and $q_{f}$ lie in a geodesically-compact region of the configuration space, the geodesics represent straightline paths between the points, and so require no lateral acceleration. The shortest geodesic is also the shortest path between the points; it thus requires the smallest average speed to traverse the curve in a given time, and therefore the least tangential acceleration to achieve this average speed while starting and ending at rest.
\item If $\configdot_{0}$ or $\configdot_{f}$ is non-zero, then the shape of the Riemannian spline solution resembles that of a ``stiff beam" whose end or ends are clamped into a fixed orientation and stretch at the endpoints. As per the definition of the cost function, the specific shape assumed by the trajectory minimizes the integrated-squared-norm of its covariant acceleration, which has the effect of making the time-progression move slowly over sections of high path-curvature (as if the ``beam" material contracts as it bends, making it denser in curved regions).

Note that this beam analogy does not extend to there being ``buckled" solutions for the trajectory when the boundary conditions require the speeds at the endpoints to be faster than the average speed over the trajectory---unlike a physical beam with a fixed length, for which buckling can relieve compressive stress, moving away from the most direct path between the specified endpoints does not ``absorb any length", and so does not reduce the amount by which the trajectory must slow down and speed up to match the boundary speeds to the total time for the trajectory.
\end{enumerate}

\subsubsection{The general cometric case: Biased Riemannian Splines and Geodesics}
\label{subsec:biasedsplines}

In this subsection, we provide a conceptual description of the biased splines and geodesics associated with the data $(\configspace,\metricg,\cometricg)$, which were theoretically characterized in Theorem~\ref{thm:solutions}.  The idea is, later, to apply these results to the case $\cometricg = \tg$ given by the torque-cost cometric of a robotics problem.

\medskip

The Hamiltonian system~\eqref{eq:bsystemh} can be re-written in a form which generalizes the force-based spline equations in~\eqref{eq:splineforce} to account for arbitrary cometrics. To this end, we replace the force $F$ covector in~\eqref{eq:splineforce} with the effort covector $\effort=h(a,\bullet)$ described in \S\ref{sec:effortcovectordescrip} (which reduces to $F$  when $\tilde g=g^*$). The second eq. in~\eqref{eq:bsystemh} then reads $\alpha=2 E$ while the third equation reads $p=-2 \nabla^*_v E$. With these, we obtain the equivalent effort-based system,%
\begin{subequations}
\label{eq:biasedspline}
\begin{align}    v &= \dot{x} &&&&&& \text{velocity} \label{eq:forcesplineforcevel}\\
a &= \nabla_{v}v && && &&  \text{covariant acceleration}\\
&& E &= \inmetricg(a,\bullet) && && \text{effort} \label{eq:forceaccelbulletdual}\\
&& && \tilde{\yank} &= \dualconnection_{v}\effort  &&  \text{effort-yank} \\
&& && \tilde{\tug} &= \dualconnection_{v}\yank =  \covprod{\effort}{-(R_g(\bullet, v) v)} + \frac{1}{2}(\dualconnection h^*)_{\bullet}(\effort,\effort). &&  \text{effort-tug} \label{eq:forcesplineforcetug}
\end{align}
\end{subequations}
The last term above can be interpreted as accounting for changes in the induced cometric $h^*$ that change the effort required to generate system force at different configurations $q\in Q$ (e.g., the leverage that the actuators have on the system components). As noted before, these changes are geometrically captured by the tensor $\bnabla h^*$ (see eq.~\eqref{eq:bnablatildeg} with $\tilde g=h^*$ there) which measures the incompatibility of the induced cometric $h^*$ with the kinetic metric geometry. Also as mentioned before, note that when $Q=\R^n$ and $g$ is Euclidean, the above system reduces to the 4th order equation~\eqref{eq:ELeqsimple} for $q(t)$ which was derived variationally in Example~\ref{ex:ELeq} (see also Remark~\ref{rmk:eqshsol} for details on the comparison to the Hamiltonian system~\eqref{eq:bsystemh}).

\smallskip

\noindent {\bf Some qualitative observations about biased splines.} For latter use, let us assume that the cometric $\tilde g=\tg$ is given as the torque cometric discussed in \S\ref{sec:torquemetric}. In this case, the instantaneous cost is calculated in the joint coordinates as
\begin{equation}
\tilde g(F,F) = \overbrace{\begin{bmatrix} \tau_{1} & \ldots & \tau_{n} \end{bmatrix}}^{F} \overbrace{\begin{bmatrix} \tau_{1} \\ \vdots \\ \tau_{n} \end{bmatrix}}^{\hphantom{\transpose{}}\transpose{F}}.
\end{equation}
The covariant acceleration and dualization operations that produce $F$ in~\eqref{eq:optcontrol}, however, remain unchanged (because they rely only on the metric and are independent of the choice of cometric). Together, these two properties of the system equations mean that the torque-optimal trajectories have the following characteristics vis-a-vis the ordinary Riemannian splines case (where $\tilde g=g^*$):
\begin{itemize}
\item If the initial and final velocities are left unspecified, then the torque-optimal trajectory between the initial and final positions is the $(Q,g)$ geodesic determined by $q_0,q_f\in Q$. This happens because acceleration is zero, $a=0$, along the geodesic, hence the force is also zero, minimizing its integral $I$.

\item If the initial and final velocities are specified, the torque-optimizing trajectories follow paths similar to the geodesics but ``biased" they are more aligned with the minor axis of $\cometric$.

In coordinates where the cometric is flat (e.g., the actuator coordinates of a mechanical system), this alignment ``blunts" the effect of a non-flat mass matrix, pulling the optimal trajectory to be closer to a straight line or simple Euclidean spline as seen in those coordinates. In coordinates where the cometric is not flat, the bias may introduce ``buckled" solutions in which it is advantageous to accept some extra length and curvature to reduce the cost of the minimum necessary accelerations dictated by the boundary conditions.
\end{itemize}

\subsection{Examples}\label{sub:exsbiased}
In this subsection, we provide several simple examples of physical systems defined by data $(\configspace,\metricg,\cometricg)$ in which we can observe and understand the departure of biased splines and geodesics from their ordinary counterparts. For these examples, we have selected two systems that minimally illustrate the effects of using cometrics that are not dual to the metric, along with a third system that demonstrates these effects in something closer to a real-world scenario:
\begin{enumerate}
\item A system with a flat metric and a non-flat cometric
\item A system with a non-flat metric and a flat cometric
\item A two-link robotic arm with two different actuation schemes.
\end{enumerate}
We present additional examples of the differences between optimizing against the dual-mass and torque-squared cometrics (considered in the context of cyclic trajectories and locomotion) in~\cite{Hatton:2022Inertial}.

\begin{example}
Our first example system is a flat manifold (Euclidean plane) equipped with a cometric for which the cost of forces in the $x$ direction increases quadratically with $y$, i.e. with $\metricg$ and $\cometricg$ instantiated as
\begin{equation}
\metric = \begin{bmatrix}
1 & 0 \\ 0 & 1
\end{bmatrix}
\qquad \text{and} \qquad
\cometric = \begin{bmatrix}
1+10y^{2} & 0 \\ 0 & 1
\end{bmatrix}
.
\end{equation}
If we specify the boundary conditions as being at two points separated by some amount in the $x$ direction and starting and ending at rest, then the acceleration-optimal trajectory under $(\metric,\dualmetric)$ is an acceleration/decelleration motion along the line connecting these points, and the torque-optimal trajectory dips towards the $y=0$ line as illustrated in Fig.~\ref{fig:quadraticcometric}.

\begin{figure}
\centering
\includegraphics[width=.8\textwidth]{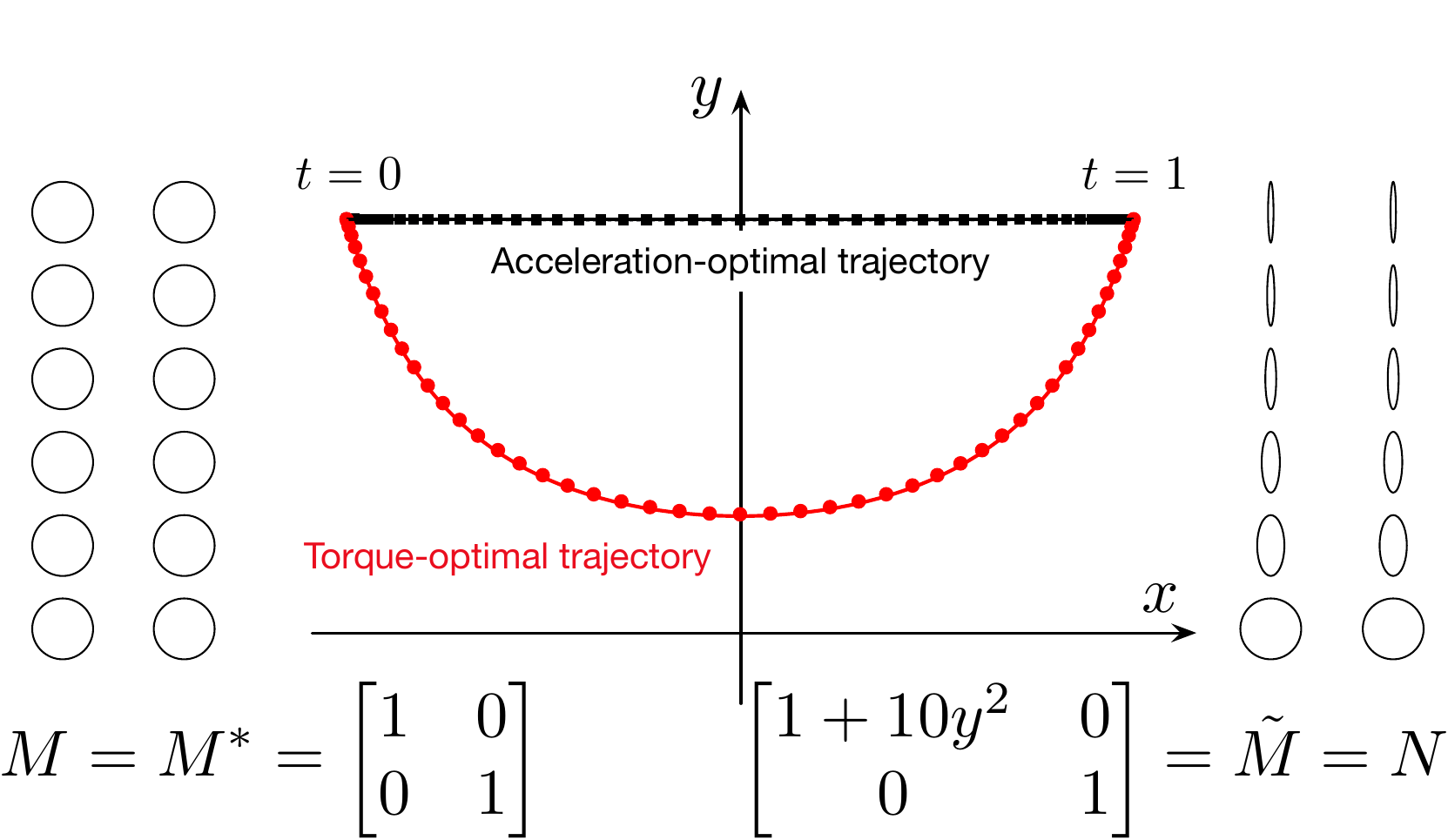}
\caption{Flat metric, quadratic cometric. On a space with a flat metric, the acceleration-optimal trajectory between two points (starting and ending at rest) is the straight-line geodesic connecting the points. If we introduce a cometric whose $xx$ component increases with $y$, then the torque optimal curve dips down, accepting a longer pathlength in exchange for a lower torque-per-acceleration cost.}
\label{fig:quadraticcometric}
\end{figure}

Note that this dip is not an effect of the optimizer trying to put the middle of the curve at a small $y$ value: there is no $x$-acceleration at the midpoint, so the cheaper $x$-acceleration at this point plays no role in the cost of the trajectory. Instead, the shape of the trajectory is a result of the ends of the trajectory (where $x$-acceleration happens) being at a lower $y$ value than the straight line, along with the property that the distortion in the cometric makes adding $y$ acceleration to an existing $x$ acceleration very inexpensive.

If, as illustrated in Fig.~\ref{fig:quadraticcometriccompressed}, we set the initial and final velocities to be parallel with the geodesic connecting the endpoints, but with speeds much greater than the geodesic speed for the given time interval, the acceleration-optimal curve suffers a hard deceleration along the geodesic path to kill time, followed by a hard acceleration to satisfy the ending boundary condition. Given sufficiently large ``compression" of the boundary conditions in this manner, the torque-optimal trajectory may ``loop", delaying  any $x$ decelleration until it can take advantage of the reduced force cost near $y=0$, and then accepting a greater total $x$-acceleration (moving with negative $x$-velocity instead of simply slowing down) in exchange for not having to apply $x$ force at large $y$.

\begin{figure}
\centering
\includegraphics[width=.8\textwidth]{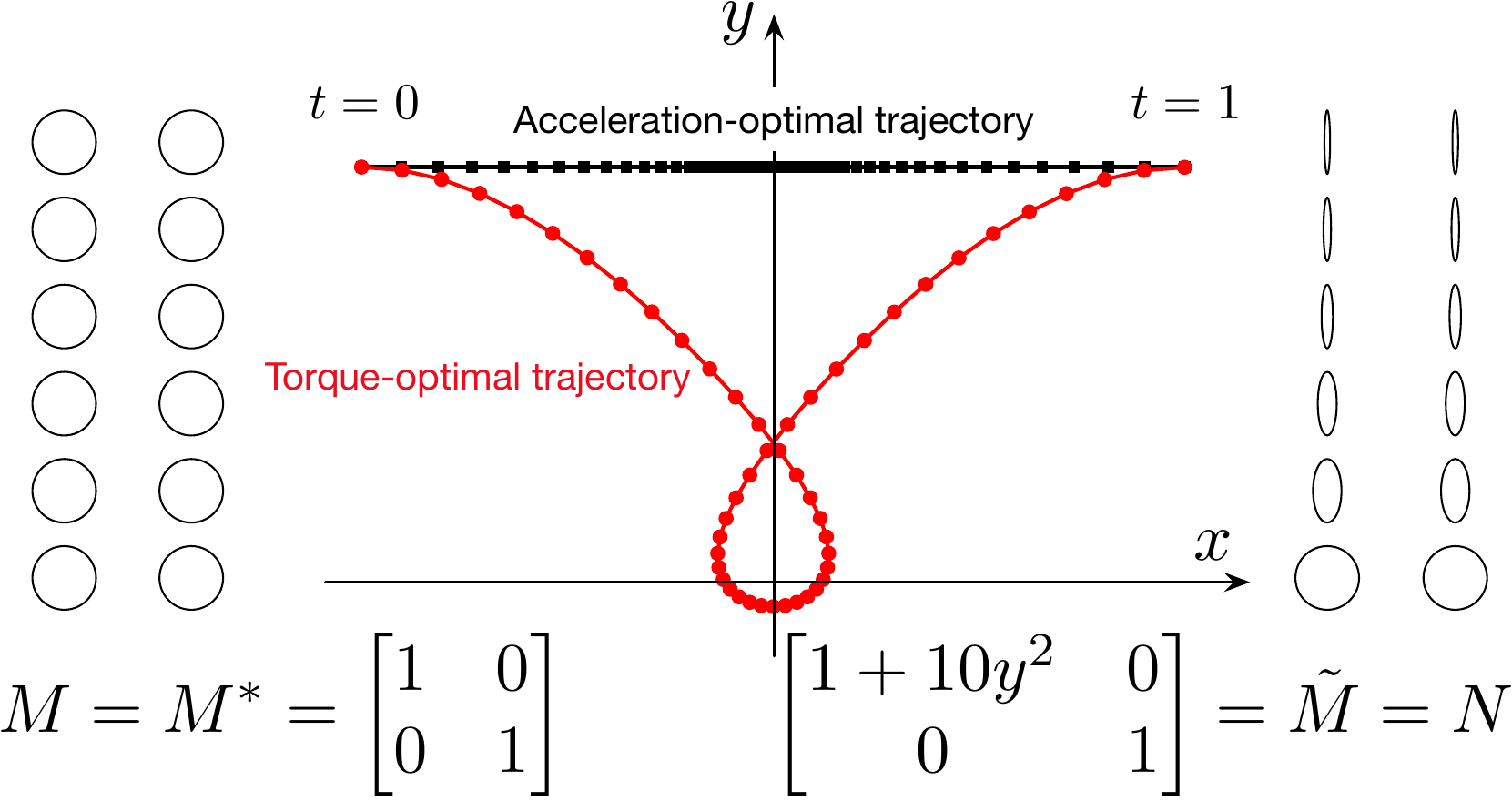}
\caption{Compressed boundary conditions with flat metric, quadratic cometric. Under highly ``compressed" boundary conditons (initial and final velocities aligned with the geodesic connecting the endpoints, but with starting and end speeds much larger than average speed), the acceleration-optimal curve slows down to kill time and then speeds back up, but the torque-optimal trajectory can deviate from this line to take advantage of cheaper acceleration costs elsewhere in the configuration space.}
\label{fig:quadraticcometriccompressed}
\end{figure}

Finally, setting the initial and final velocities to leave and return to the geodesic between the endpoints, as illustrated in Fig.~\ref{fig:quadraticcometricscurve}, demonstrates that the offset from the geodesic spline to the biased spline can be in different directions based on the structure of the curve: When the geodesic spline is entirely in a positive-$y$ region, the biased spline is shifted ``down" towards $y=0$, but when the geoedsic spline crosses the $y=0$ line, the bias is ``inside" the geodesic spline rather than a bulk offset.

\begin{figure}
\centering
\includegraphics[width=.8\textwidth]{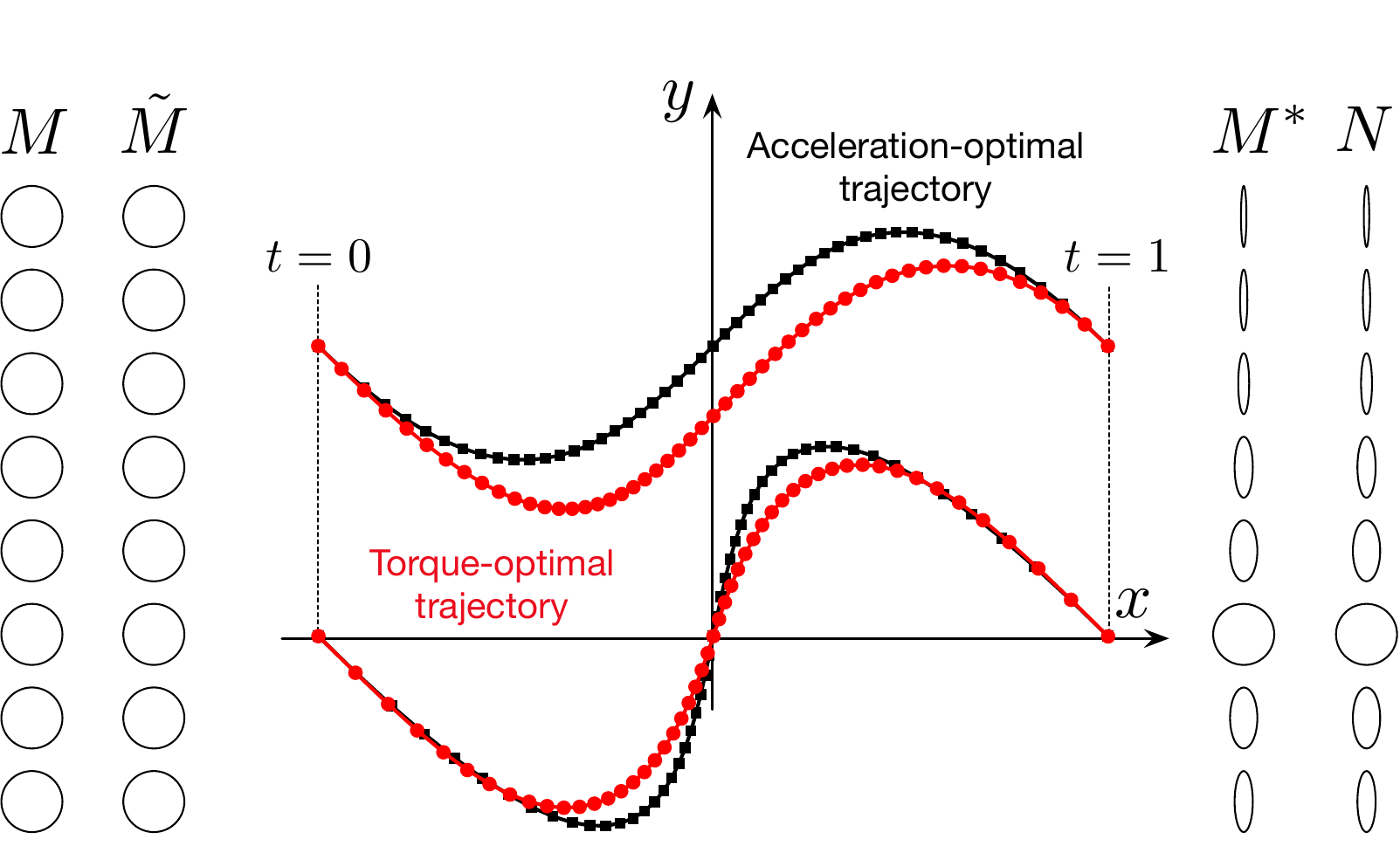}
\caption{Skewed boundary conditions with flat metric, quadratic cometric. If the initial and final velocity are not aligned with the geodesic between the endpoints, the bias in the optimal torque trajectories towards low-cost regions may be in different directions at different points along the curve.}
\label{fig:quadraticcometricscurve}
\end{figure}

\end{example}

\begin{example}
Our second example replaces the dual metric on the sphere with a cometric that is flat in longitude-latitude coordinates (corresponding to having a yaw-pitch mechanism with equal motor costs). As illustrated in Fig.~\ref{fig:spheremetricflatcometric}, this change results in torque-optimal trajectories that rise less in latitude than do the great-circle geodesics between the endpoints.

\begin{figure}
\centering
\includegraphics[width=\textwidth]{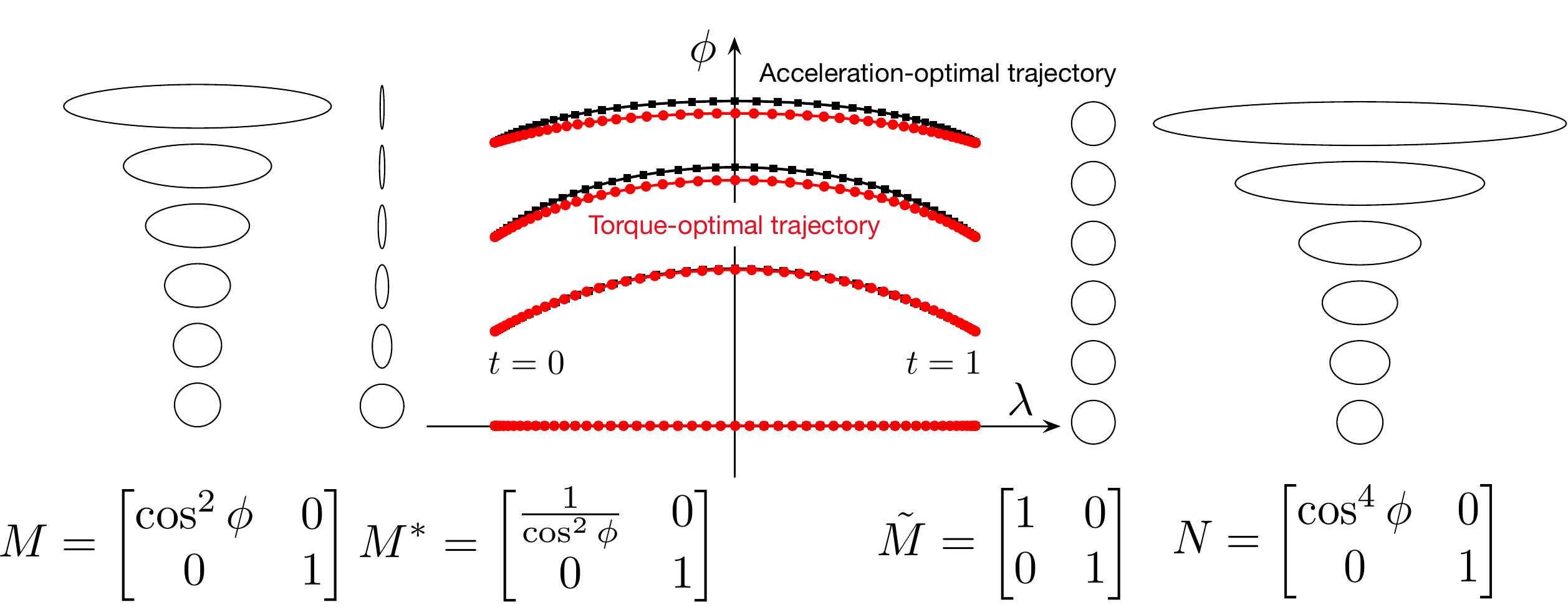}
\caption{Sphere metric, flat cometric. Flattening the cometric on the sphere removes some of the incentive to add latidudinal acceleration into required longitudinal acceleration, so the torque-optimal trajectories are flatter than the geodesic-following acceleration-optimal trajectories.}
\label{fig:spheremetricflatcometric}
\end{figure}

Given that the change from $\dualmetric$ to $\cometric$ makes $\lambda$-acceleration less expensive at large $\phi$, it might be intuitive to expect that the torque-optimal curves would rise to larger latitude than the geodesics, rather than smaller latitudes. Instead, however, the change in the trajectory can be better understood by noting that the horizontal stretch of the $\cometric$ indicatrices relative to the $\dualmetric$ indicatrices (and of the $\inmetric$ indicatrices relative to those of $\metric$) means that at any point on the geodesic spline that is accelerating in both the $\lambda$ and $\phi$ directions, the biased spline can afford to increase its $\lambda$ acceleration, which flattens the trajectory. We also note that being at a large $\phi$ in the middle of the trajectory is not directly useful in reducing the torque-cost of the trajectory, because there is no $\lambda$-torque cost at that point.

\end{example}

\begin{example}

Finally, we consider the motion of a two-link arm whose mass is concentrated at the distal end of the link. As illustrated in Fig.~\ref{fig:serialparallelspline}, the acceleration-optimal trajectories for moving this point-mass follow straight lines through the world. Imposing the actuator-torque cometrics for serial actuation (a motor controlling the relative angle at each joint) or parallel actuation (a motor controlling the absolute orientation of each link) both result in torque-optimal curves that are biased away from these straight lines, with the bias more pronounced in the parallel-actuated mechanism, as would be expected from the larger aspect ratio in the world-space indicatrices for $\cometric$ and $\inmetric$ illustrated in the middle-bottom plots in Fig.~\ref{fig:twolink}.

\begin{figure}
\centering
\includegraphics[width=.8\textwidth]{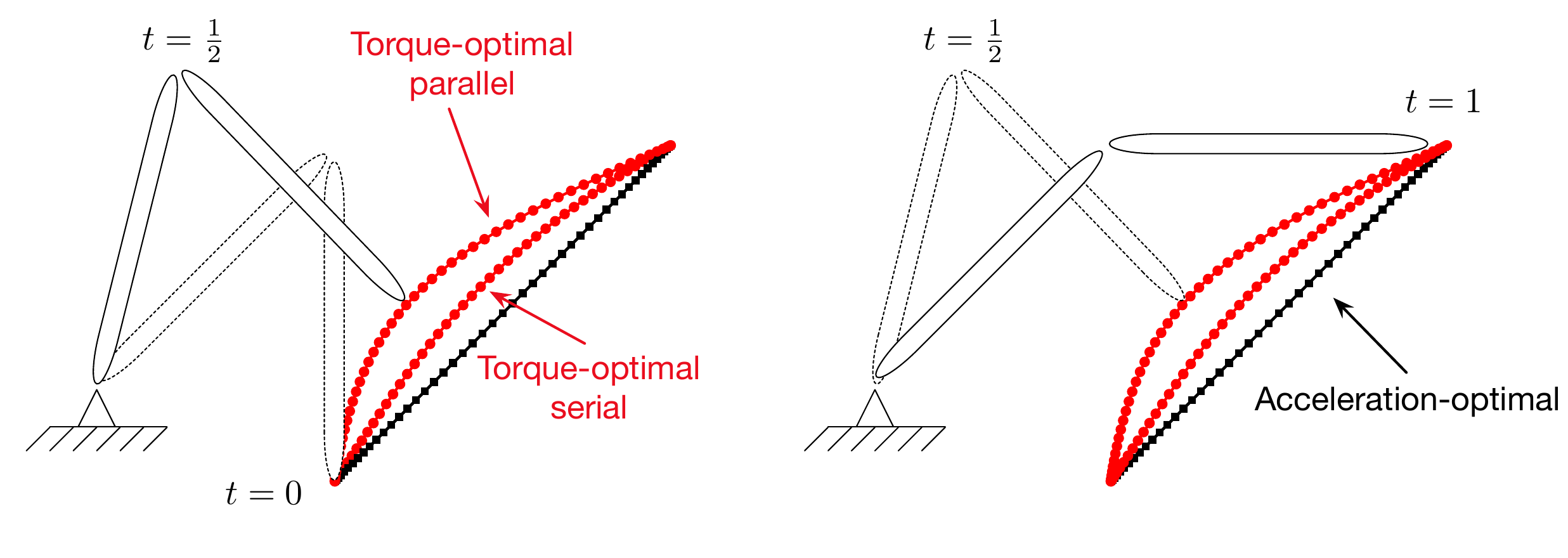}
\caption{Serial- and parallel-actuated two-link arms. The acceleration-optimal trajectory for a two-link arm with mass concentrated at the end of the arm is to move the mass in a straight line. Using the cometrics corresponding to torques on the actuators leads to curves that initially move the mass more in the $y$ direction and then in the $x$ direction, with this effect more pronounced for parallel actuation.}
\label{fig:serialparallelspline}
\end{figure}
\end{example}

\appendix
\section{Some coordinate formulas} \label{app:coordinates}

Let $U\subset \R^\qdim$ be an open subset with coordinates $q=(q^1,...,q^\qdim)$. We shall think of $U$ as a coordinate domain in a configuration manifold $\configspace$ and write down explicit coordinate formulas for some of the tensors that appear in our mechanical systems.

First, $(\partial_{q^i})$ denotes the basis of $TU$ given by (constant) vector fields in $U$ and $(dq^i)$ the dual basis of $T^*U$. A metric $\metricg$ on $U$ is thus defined by its coeffient functions
\[ g_{ij}(q) = \metricg(\partial_{q^i},\partial_{q^j})|_q .\]
Similarly, a cometric $\cometricg$ on $U$ we have its coefficient functions
\[ \cometricg^{ij}(q) = \cometricg(dq^i,dq^j)|_q.\]
When $\cometricg = \dualmetricg$ is the dual of $\metricg$, we have that the matrix $\cometric = (\cometricg^{ij}(q))$ is the inverse of $M=(g_{ij}(q))$ at each $q\in U$.

Let $\nabla$ be the Levi-Civitta connection on $TU$ induced by $\metricg$ (or $\nabla$ any linear connection on $TU$). Its Christoffel symbols $\Gamma_{ij}^k\equiv \Gamma_{ij}^k(q)$ are defined by
\[ \nabla_{\partial_{q^i}} \partial_{q^j} = \Gamma_{ij}^k \partial_{q^k},\]
where Einstein's summation convention will be used throughout. Note that there are well-known explicit formulas for the $\Gamma_{ij}^k$ in terms of the $g_{ij}$'s and their derivatives when $\nabla$ is the Levi-Civitta connection:
\[ g_{ij}\Gamma^i_{kl} = \frac{1}{2} \left( \partial_{q^l} g_{jk} + \partial_{q^k}g_{jl} - \partial_{q^j}g_{kl} \right). \]
The dual connection $\dualconnection$ on $T^*U$  satisfies
\[ \dualconnection_{\partial_{q^i}} dq^j = -\Gamma_{ik}^j dq^k.\]
This is because $\dualconnection$ is defined by the relation $L_X \langle \alpha, Y\rangle = \langle \dualconnection_X \alpha, Y\rangle + \langle \alpha, \nabla_X Y\rangle$, for all vector fields $X,Y$ and all 1-forms $\alpha$ (and we can take $X=\partial_{q^i}, \ Y=\partial_{q^j} \ \alpha = dx^k$).

The Riemann curvature $R^\nabla\in \Omega^2(U, End(TU))$ will be described by its tensorial coefficient functions $R_{ijk}^l\equiv R_{ijk}^l(q)$ defined by
\[\langle \alpha, R^\nabla(X,Y)Z \rangle = R_{ijk}^l X^i Y^j Z^k \alpha_l, \]
where $X=X^i\partial_{q^i}, \ Y= Y^j \partial_{q^j}, \ Z=Z^k\partial_{q^k}$ are any vectors and $\alpha = \alpha_l dq^l$ any covector. There is an explicit well-known formula for the $R_{ijk}^l$ in terms of the $\Gamma_{ij}^k$'s that follows from the definition of $R^\nabla$:
\[ R_{ijk}^l = \partial_{q^i} \Gamma^l_{jk} - \partial_{q^j}\Gamma^l_{ik} + \Gamma^n_{jk} \Gamma^l_{in} - \Gamma^n_{ik} \Gamma^l_{jn}. \]
(This formula can be specialized further when $\nabla$ is the Levi-Civitta connection, but this expression will not be needed here.)

\smallskip

With the above conventions, the \ctensor\ $\dualconnection \cometricg$ is then given by
\[ (\dualconnection \cometricg)_X (\alpha,\beta)|_q = \tau^{ij}_k(q) X^k \alpha_i \beta_j, \text{ for any } X=X^k \partial_{q^k}, \ \alpha = \alpha_i dq^i, \ \beta_j dq^j,\]
where the relevant coeffients $\tau^{ij}_k(q)$ are given by
\[ \tau^{jk}_i(q) = \partial_{q^i}\cometricg^{jk} + \Gamma_{il}^j\cometricg^{lk} + \Gamma_{il}^k \cometricg^{jl}. \]
This follows directly from the definition of $\dualconnection \cometricg$ in eq.~\eqref{eq:bnablatildeg} (taking $X=\partial_{q^k}$, $\omega_1=dx^i$ and $\omega_2=dx^j$ there).

\medskip

Finally, let us now write down in coordinates the system of Hamiltonian equations~\eqref{eq:bsystemh} for the biased splines. The unknowns are: $\gamma(t)\equiv q(t) \in U$ a curve in the coordinate domain, $v(t)=\configdot(t) \in T_{q(t)}U = \R^\qdim$ its velocity, and two 1-forms (along $q(t)$) $\alpha(t) = \alpha_i(t) dq^i$ and $p(t) = p_i(t) dq^i$. Note that the true unkowns behind $\alpha(t)$ and $p(t)$ are just the real-valued coeffiecients $(\alpha_i(t))$ and $(p_i(t))$ (the condition of the 1-forms being defined along $q(t)$ poses no constraint in the coordinate situation).
With these notations, the eqs. in system~\eqref{eq:bsystemh} are then equivalent to the following system of scalar equations indexed by a given $i=1,...,\qdim$ (all the other indices appearing are summation ones),
\begin{eqnarray}
\configdot^i &=& v^i \nonumber \\
\dot v^i + \Gamma^i_{jk}\configdot^j v^k &=& \frac{1}{2} (h^*)^{ji}\alpha_j \nonumber \\
\dot \alpha_i - \Gamma^j_{ki} \configdot^k \alpha_j &=& - p_i \nonumber \\
\dot p_i - \Gamma^j_{ki} \configdot^k p_j &=& -\frac{1}{4} \tau_i^{jk} \alpha_j \alpha_k + R_{ijk}^l v^j v^k \alpha_l, \label{eq:coordhsystem}
\end{eqnarray}
where $h^*$ is the cometric dual to the metric $h$ defined in~\eqref{eq:defhfromtg} and $\tau_i^{jk}$ are the coefficients of $\dualconnection h^*$ defined above (with $\tilde g=h^*$).

\bibliographystyle{siamplain}

\end{document}